\documentclass{amsart}
\usepackage{amssymb,amscd}
\usepackage{color}
\usepackage{enumerate}

\usepackage[arrow,matrix]{xy}
\newcommand\datver[1]{\def\datverp%
 {\par\boxed{\boxed{\text{Version: #1; Run: \today}}}}}
\datver{3.2; Revised: 18/07/2005}

\newcommand{\End}{\operatorname{End}}

\newcommand{\pa}{\partial}

\newcommand{\supp}{\operatorname{supp}}

\newcommand{\CI}{\mathcal{C}^\infty}

\newcommand{\bS}{{}^b\kern-1pt S}
\newcommand{\bT}{{}^b\kern-1pt T}
\newcommand{\Hom}{\operatorname{Hom}}

\newcommand{\cC}{\mathcal{C}}

\newcommand{\cF}{\mathcal{F}}

\newcommand{\cmun}{\cC^{-\infty}}
\newcommand{\codim}{\operatorname{depth}}

\newcommand{\cR}{\mathcal{R}}
\newcommand{\cun}{\cC^{\infty}}

\newcommand{\cunc}{\cun_{c}}

\newcommand{\ffi}{\cF_{{\rm fiber}}^{-1}}

\newcommand{\nzn}{\NN_{0}}
\newcommand{\oeh}{\Omega^{1/2}}

\newcommand{\piM}{\pi_{M}}

\newcommand{\sia}[1]{\sigma^{(#1)}}

\newcommand{\CC}{\mathbb C}
\newcommand{\NN}{\mathbb N}
\newcommand{\RR}{\mathbb R}

\newcommand{\ZZ}{\mathbb Z}

\newcommand{\CIc}{{\mathcal C}^{\infty}_{\text{c}}}

\newcommand{\GR}{\mathcal G}

\newcommand{\be }{\begin{eqnarray*}}
\newcommand{\ee }{\end{eqnarray*}}

\newcommand{\mF}{\mathcal F}

\newcommand{\mR}{\mathcal R}

\newcommand{\Lie}{\operatorname{Lie}}
\newcommand{\fgg}{\mathfrak g}

\newcommand\comment[1]{}

\let\<\langle
\let\>\rangle

\let\rho\varrho

\newcommand{\ie}{{\em i.e., }}
\newcommand{\VV}{\mathcal V}
\newcommand{\Diff}[1]{{\rm Diff}(#1)}
\newcommand{\DiffV}[1]{{\rm Diff}^{#1}_{\VV}}
\newtheorem{theorem}{Theorem}[section]
\newtheorem{proposition}[theorem]{Proposition}
\newtheorem{corollary}[theorem]{Corollary}

\newtheorem{lemma}[theorem]{Lemma}

\theoremstyle{definition}
\newtheorem{definition}[theorem]{Definition}
\theoremstyle{remark}
\newtheorem{remark}[theorem]{Remark}

\newtheorem{examples}[theorem]{Examples}

\def\supp{\text{supp} }

\long\def\komment#1{}

\def\be{{\beta}}

\def\phi{{\varphi}}

\def\pa{{\partial}}

\def\vol{{\mathop{{\rm vol}}}}

\author[B. Ammann]{Bernd Ammann} \address{Bernd Ammann, 
L'Institut \'Elie Cartan (I\'ECN),
Universit\'e Henri Poincar\'e, B.P. 239, F-54506 Vandoeuvre-Les-Nancy, Cedex,
France
}
\email{ammann@iecn.u-nancy.fr}

\author[R. Lauter]{Robert Lauter} \address{Robert Lauter, Universit\"at Mainz,
Fachbereich 17 -- Mathematik, D - 55099 Mainz, Germany}
\email{lauter@mathematik.uni-mainz.de}

\author[V. Nistor]{Victor Nistor} \address{Victor Nistor, Pennsylvania State
       University, Math. Dept., University Park, PA 16802}
       \email{nistor@math.psu.edu}

\thanks{Ammann was partially supported by the European Contract Human 
Potential Program, Research Training Networks HPRN-CT-2000-00101 
and HPRN-CT-1999-00118,
Nistor was partially supported by the NSF Grants DMS-9971951
and DMS-0200808. Manuscripts available from {\bf
http:{\scriptsize//}www.math.psu.edu{\scriptsize/}nistor{\scriptsize/}}
.}

\begin{document}

\dedicatory\datverp

\title[Pseudodifferential operators]{Pseudodifferential operators on
manifolds with a Lie structure at infinity}

\begin{abstract}
We define and study an algebra $\Psi_{1,0,\VV}^\infty(M_0)$ of
pseudodifferential operators canonically associated to a non-compact,
Riemannian manifold $M_0$ whose geometry at infinity is described by a
Lie algebra of vector fields $\VV$ on a compactification $M$ of $M_0$
to a compact manifold with corners.  We show that the basic properties
of the usual algebra of pseudodifferential operators on a compact
manifold extend to $\Psi_{1,0,\VV}^\infty(M_0)$. We also consider the
algebra $\DiffV{*}(M_0)$ of differential operators on $M_0$ generated
by $\VV$ and $\CI(M)$, and show that $\Psi_{1,0,\VV}^\infty(M_0)$ is a
microlocalization of $\DiffV{*}(M_0)$. Our construction solves a
problem posed by Melrose in 1990. Finally, we introduce and study
semi-classical and ``suspended'' versions of the algebra
$\Psi_{1,0,\VV}^\infty(M_0)$.
\end{abstract}

\maketitle \tableofcontents

\section*{Introduction}

Let $(M_0,g_0)$ be a complete, non-compact Riemannian manifold.  It is a
fundamental problem to study the geometric operators on $M_0$.  As in the
compact case, pseudodifferential operators provide a powerful tool for that
purpose, provided that the geometry at infinity is taken into account.  One
needs, however, to restrict to suitable classes of non-compact manifolds.

Let $M$ be a compact manifold with corners such that $M_0 = M \smallsetminus
\pa M$, and assume that the geometry at infinity of $M_0$ is described by a
Lie algebra of vector fields $\VV \subset \Gamma(M; TM)$, that is, $M_0$ is a
Riemannian manifold with a Lie structure at infinity,
Definition~\ref{def.R.usi}.  In \cite{meicm}, Melrose has formulated a far
reaching program to study the analytic properties of geometric differential
operators on $M_0$. An important ingredient in Melrose's program is to define
a suitable pseudodifferential calculus $\Psi_\VV^\infty(M_0)$ on $M_0$ adapted
in a certain sense to $(M, \VV)$. This pseudodifferential calculus was called
a ``microlocalization of $\DiffV{*}(M_0)$'' in \cite{meicm}, where
$\DiffV{*}(M_0)$ is the algebra of differential operators on $M_0$ generated
by $\VV$ and $\CI(M)$. (See Section \ref{sec.KN}.)

Melrose and his collaborators have constructed the algebras
$\Psi_\VV^\infty(M_0)$ in many special cases, see for instance \cite{emm91,
  Mazzeo, mame87, MaMeAsian, me81, meaps, MelroseMendoza, jaredduke}, and
especially \cite{MelroseScattering}.  One of the main reasons for considering
the compactification $M$\ is that the geometric operators on manifolds with a
Lie structure at infinity identify with degenerate differential operators
on~$M$.  This type of differential operators appear naturally, for example,
also in the study of boundary value problems on manifolds with singularities.
Numerous important results in this direction were obtained also by Schulze and
his collaborators, who typically worked in the framework of the Boutet de
Monvel algebras. See \cite{ ScSc,schwil} and the references therein.  Other
important cases in which this program was completed can be found in \cite{zfr,
  defr, LN1, nwx, Parenti}.  An earlier important motivation for the
construction of these algebras was the method of layer potentials for boundary
value problems and questions in analysis on locally symmetric spaces. See for
example \cite{BNXZ, Baer, Costabel, Dauge, Lewis, LewisParenti, MazzeoVasy,
  MiNi}.

An outline of the construction of the algebras $\Psi_\VV^\infty(M_0)$ was
given by Melrose in \cite{meicm}, provided certain compact manifolds with
corners (blow-ups of $M^2$ and $M^3$) can be constructed.  In the present
paper, we modify the blow-up construction using Lie groupoids, thus completing
the construction of the algebras $\Psi_\VV^\infty(M_0)$.
Our method relies on recent progress achieved in 
\cite{aln1, CrainicFernandes, nwx}.

The explicit construction of the algebra $\Psi_{1,0,\VV}^{\infty}(M_0)$
microlocalizing $\DiffV{*}(M_0)$ in the sense of \cite{meicm} is, roughly, as
follows. First, $\VV$ defines an extension of $TM_0$ to a vector bundle $A \to
M$ ($M_0 = M \smallsetminus \pa M$). Let $V_r := \{d(x,y) < r\} \subset M_0^2$
and $(A)_r = \{v \in A, \|v\| < r\}$. Let $r>0$ be less than the injectivity
radius of $M_0$ and $V_r \ni (x, y) \mapsto (x, \tau(x, y)) \in (A)_r$ be a
local inverse of the Riemannian exponential map $TM_0 \ni v \mapsto \exp_x(-v)
\in M_0 \times M_0$. Let $\chi$ be a smooth function on $A$ with support in
$(A)_r$ and $\chi = 1$ on $(A)_{r/2}$. For any $a \in S^m_{1,0}(A^*)$, we
define
\begin{equation}
  \big[a_i(D)u\big](x) = (2\pi )^{-n} \int_{M_{0}} \left
    (\int_{T^{*}_{x}M_{0}} e^{i \tau(x,y) \cdot \eta} \chi(x,
    \tau(x,y)) a(x, \eta)u(y)\, d \eta \right) dy.
\end{equation}
The algebra $\Psi_{1,0,\VV}^{\infty}(M_0)$ is then defined as the linear span
of the operators $a_\chi(D)$ and $b_\chi(D)\exp(X_1) \ldots \exp(X_k)$, $a \in
S^\infty(A^*)$, $b \in S^{-\infty}(A^*)$, and $X_j \in \VV$, and where
$\exp(X_j) : \CIc(M_0) \to \CIc(M_0)$ is defined as the action on functions
associated to the flow of the vector field $X_j$.
 
The operators $b_\chi(D)\exp(X_1)\ldots \exp(X_k)$ are needed to make our
space closed under composition.  The introduction of these operators is in
fact a crucial ingredient in our approach to Melrose's program. The results of
\cite{CrainicFernandes, nwx} are used to show that
$\Psi_{1,0,\VV}^{\infty}(M_0)$ is closed under composition, which is the most
difficult step in the proof.

A closely related situation is encountered when one considers a product of a
manifold with a Lie structure at infinity $M_0$ by a Lie group $G$ and
operators $G$ invariant on $M_0 \times G$. We obtain in this way an algebra
$\Psi_{1,0,\VV}^{\infty}(M_0; G)$ of $G$--invariant pseudodifferential
operators on $M_0 \times G$ with similar properties. The algebra
$\Psi_{1,0,\VV}^{\infty}(M_0; G)$ arises in the study of the analytic
properties of differential geometric operators on some higher dimensional
manifolds with a Lie structure at infinity.  When $G = \RR^q$, this algebra is
slightly smaller than one of Melrose's suspended algebras and plays the same
role, namely, it appears as a quotient of an algebra of the form
$\Psi_{1,0,\VV'}^{\infty}(M_0')$, for a suitable manifold~$M_0'$. The quotient
map $\Psi_{1,0,\VV'}^{\infty}(M_0') \to \Psi_{1,0,\VV}^{\infty}(M_0; G)$ is a
generalization of Melrose's indicial map. A convenient approach to indicial
maps is provided by groupoids \cite{LN1}.

We also introduce a semi-classical variant of the algebra
$\Psi_{1,0,\VV}^{\infty}(M_0)$, denoted $\Psi_{1,0,\VV}^{\infty}(M_0[[h]])$,
consisting of semi-classical families of operators in
$\Psi_{1,0,\VV}^{\infty}(M_0)$.  For all these algebras we establish the usual
mapping properties between appropriate Sobolev spaces.

The article is organized as follows. In Section \ref{sec.LI} we recall the
definition of manifolds with a Lie structure at infinity and some of their
basic properties, including a discussion of compatible Riemannian metrics. In
Section \ref{sec.KN} we define the spaces $\Psi_{1,0,\VV}^{m}(M_0)$ and the
principal symbol maps.  Section \ref{sec.product} contains the proof of the
crucial fact that $\Psi_{1,0,\VV}^{\infty}(M_0)$ is closed under composition,
and therefore it is an algebra. We do this by showing that
$\Psi_{1,0,\VV}^{\infty}(M_0)$ is the homomorphic image of $\Psi_{1,
  0}^\infty(\GR)$, where $\GR$ is any $d$-connected Lie groupoid integrating
$A$ ($d$--connected means that the fibers of the domain map $d$ are
connected).  In Section \ref{sec.Prop} we establish several other properties
of the algebra $\Psi_{1,0,\VV}^{\infty}(M_0)$ that are similar and analogous
to the properties of the algebra of pseudodifferential operators on a compact
manifold.  In Section \ref{Sec.Ginv} we define the algebras
$\Psi_{1,0,\VV}^{\infty}(M_0[[h]])$ and $\Psi_{1,0,\VV}^{\infty}(M_0;G)$,
which are generalizations of the algebra $\Psi_{1,0,\VV}^{\infty}(M_0)$. The
first of these two algebras consists of the semi-classical (or adiabatic)
families of operators in $\Psi_{1,0,\VV}^{\infty}(M_0)$. The second algebra is
a subalgebra of the algebra of $G$--invariant, properly supported
pseudodifferential operators on $M_0 \times G$, where $G$ is a Lie group.

\subsubsection*{Acknowledgements:}\ We thank Andras Vasy for several
interesting discussions and for several contributions to this paper.  R. L. is
grateful to Richard B.\ Melrose for numerous stimulating conversations and
explanations on pseudodifferential calculi on special examples of manifolds
with a Lie structure an infinity.  V. N. would like to thank the Institute
Erwin Schr\"odinger in Vienna and University Henri Poincar\'e in Nancy, where
parts of this work were completed.


\section{Manifolds with a Lie structure at infinity}
\label{sec.LI}

For the convenience of the reader, let us recall the definition of a
Riemannian manifold with a Lie structure at infinity and some of its
basic properties.

\subsection{Preliminaries}\
In the sequel, by a {\em manifold} we shall always understand a
$C^\infty$-manifold {\em possibly with corners}, whereas a {\em smooth
manifold} is a $C^\infty$-manifold {\em without corners} (and without
boundary). By definition, every point $p$ in a manifold with corners $M$ 
has a coordinate neighborhood diffeomorphic to $[0,\infty)^k \times
\RR^{n-k}$ such that the transition functions are smooth up to the
boundary. If $p$ is mapped by this diffeomorphism to $(0, \ldots, 0, 
x_{k+1}, \ldots, x_n)$, we shall say that $p$ is
a point of {\em boundary depth
$k$} and write $\codim(p) = k$. The closure of a connected component of
points of boundary depth $k$ is called a {\em face of codimension
$k$}.  Faces of codimension~$1$ are also called {\em hyperfaces}.  For
simplicity, we always assume that each hyperface $H$ of a manifold
with corners $M$ is an embedded submanifold and has a defining
function, that is, that there exists a smooth function $x_H \ge 0$ on
$M$ such that $$ H = \{ x_H = 0 \} \; \text{ and }\; dx_H \not = 0
\;\text{ on }\; H.
$$ For the basic facts on the analysis on manifolds with corners we
refer to the forthcoming book \cite{MelroseCorners}. We shall denote by $\pa
M$ the union of all non-trivial faces of $M$ and by $M_0$
the interior of $M$, \ie\ $M_0:= M \smallsetminus \pa M$.
Recall that a map $f : M \to N$ is a {\em submersion of manifolds with
corners} if $df$ is surjective at any point and $df_p(v)$ is an 
inward pointing vector if, and only if, $v$ is an inward pointing
vector. In particular, the sets $f^{-1}(q)$ are smooth manifolds (no boundary
or corners).

To fix notation, we shall denote the sections of a vector
bundle $V \to X$ by $\Gamma(X, V)$, unless $X$ is understood, in
which case we shall write simply $\Gamma(V)$. 
A Lie subalgebra $\VV\subseteq \Gamma(M, TM)$ of the Lie algebra
of all smooth vector fields on $M$ is said to be {\em a structural
Lie algebra of vector fields} provided it is a finitely generated,
projective $\cun(M)$-module and each $V\in\VV $ is tangent to all
hyperfaces of $M$.

\begin{definition}\label{def.unif.str}\
A {\em Lie structure at infinity} on a smooth manifold $M_0$ is a
pair $(M,\VV)$, where $M$ is a compact manifold, possibly with
corners, and $\VV \subset \Gamma(M, TM)$ is a
structural Lie algebra of vector fields on $M$ with the following
properties:
\begin{enumerate}[(a)]
\item $M_0$ is diffeomorphic to the interior $M \smallsetminus \pa M$
of $M,$
\item for any vector field $X$ on $M_0$ and any $p\in M_0$, 
there is a neighborhood $V$ of $p$ in $M_0$ and a vector field $Y\in \VV$,
such that $Y=X$ on $V$.
\end{enumerate}
A manifold with a  Lie structure at infinity will also be called a 
\emph{Lie manifold}.
\end{definition}

Here are some examples.
\begin{examples}\label{ex1}\ \\[-5mm] 
\begin{enumerate}[(a)] 
\item\ Take $\VV_b$ to be the set of all vector fields tangent to 
all faces of a manifold with corners $M$. Then $(M, \VV_b)$ is a 
manifold with a Lie structure at infinity.
\item\ Take $\VV_0$ to be the set of all vector fields vanishing on 
all faces of a manifold with corners $M$. Then $(M, \VV_0)$ is a 
Lie manifold. If $\pa M$ is a smooth manifold (\ie if $M$ 
is a manifold with boundary), then $\VV_0 = r\Gamma(M; 
TM)$, where  $r$ is the distance to the boundary. 
\item\ As another example consider a manifold with smooth boundary and 
consider the vector fields $\VV_{\rm sc}=r \VV_b$, where $r$ and $\VV_b$ are 
as in the previous examples.
\end{enumerate} 
\end{examples}

These three examples are, respectively, the ``$b$-calculus'', the
``$0$-calculus,'' and the ``scattering calculus'' from
\cite{MelroseScattering}. These examples are typical and will be referred to
again below.  Some interesting and highly non-trivial examples of Lie
structures at infinity on $\RR^n$ are obtained from the $N$-body problem
\cite{VasyN} and from strictly pseudoconvex domains \cite{mendoza}. Further
examples of Lie structures at infinity were discussed in \cite{aln1}.

If $M_0$ is compact without boundary, then it follows from the above
definition that $M = M_0$ and $\VV = \Gamma(M,TM)$, so a Lie structure at
infinity on $M_0$ gives no additional information on $M_0$. The interesting
cases are thus the ones when $M_0$ is non-compact.

Elements in the enveloping algebra $\DiffV{*}(M)$ of $\VV$ are called {\em
  $\VV$-differential operators on $M$}.  The order of differential operators
induces a filtration $\DiffV{m}(M)$, $m\in\nzn$, on the algebra
$\DiffV{*}(M)$. Since $\DiffV{*}(M)$ is a $\cun(M)$-module, we can introduce
$\VV$-differential operators acting between sections of smooth vector bundles
$E, F \rightarrow M$, $E, F \subset M \times \CC^N$ by
\begin{equation}\label{vectorvalued}
    \DiffV{*}(M;E,F) := e_F M_N(\DiffV{*}(M)) e_E\,,
\end{equation}
where $e_E, e_F \in M_N(\CI(M))$ are the projections onto $E$ and,
respectively, $F$. It follows that $\DiffV{*}(M;E,E) = : \DiffV{*}(M;E)$ is an
algebra that is closed under adjoints.

Let $A\rightarrow M$ be a vector bundle and $\varrho: A \rightarrow TM$ a
vector bundle map. We shall also denote by $\varrho$ the induced map
$\Gamma(M,A)\to \Gamma(M,TM)$ between the smooth section of these bundles.
Suppose a Lie algebra structure on $\Gamma(M,A)$ is given.  Then the pair
$(A,\varrho)$ together with this Lie algebra structure on $\Gamma(A)$ is
called a {\em Lie algebroid} if $\varrho([X,Y]) = [\varrho(X), \varrho(Y)]$
and $[X, fY] = f[X,Y] + (\varrho(X) f)Y$ for any smooth sections $X$ and $Y$
of $A$ and any smooth function $f$ on $M$. The map $\varrho : A \to TM$ is
called the {\em anchor of} $A$. We have also denoted by $\varrho$ the induced
map $\Gamma(M, A) \rightarrow \Gamma(M, TM)$.  We shall also write $Xf
:=\varrho(X) f$.

If $\VV$ is a structural Lie algebra of vector fields, then $\VV$ is
projective, and hence the Serre-Swan theorem \cite{Karoubi} shows that there
exists a smooth vector bundle $A_\VV\rightarrow M$ together with a natural map
\begin{equation}\label{eq.anchor}
    \begin{array}{rcccl}
      \varrho_{\VV}: A_\VV \kern-3mm& &\longrightarrow& &\kern-3mm
      TM \\ &\searrow\kern-3mm&&\kern-3mm \swarrow& \\ &&M
    \end{array}
\end{equation}
such that $\VV = \varrho(\Gamma(M, A_\VV))$.  The vector bundle $A_{\VV}$
turns out to be a {\em Lie algebroid} over $M$.

We thus see that there exists an equivalence between structural Lie algebras
of vector fields $\VV = \Gamma(A_{\VV})$ and Lie algebroids $\varrho : A
\rightarrow TM$ such that the induced map $\Gamma(M, A) \to \Gamma(M, TM)$ is
injective and has range in the Lie algebra $\VV_b(M)$ of all vector fields
that are tangent to all hyperfaces of $M$.  Because $A$ and $\VV$ determine
each other up to isomorphism, we sometimes specify a Lie structure at infinity
on $M_0$ by the pair $(M,A)$.  The definition of a manifold with a Lie
structure at infinity allows us to identify $M_0$ with $M \smallsetminus \pa
M$ and $A\vert_{M_0}$ with $TM_0$.

We now turn our attention to Riemannian structures on $M_0$.  Any metric on
$A$ induces a metric on $TM_0 = A\vert_{M_0}$.  This suggests the following
definition.

\begin{definition}
  \label{def.R.usi}\ A manifold $M_0$ with a Lie structure at infinity
  $(M,\VV)$, $\VV=\Gamma(M,A)$, and with metric $g_{0}$ on $TM_0$ obtained
  from the restriction of a metric $g$ on $A$ is called a {\em Riemannian
    manifold with a Lie structure at infinity}.
\end{definition}

The geometry of a Riemannian manifold $(M_{0},g_{0})$ with a Lie structure
$(M,\VV)$ at infinity has been studied in \cite{aln1}.  For instance,
$(M_{0},g_{0})$ is necessarily of infinite volume and complete.  Moreover, all
the covariant derivatives of the Riemannian curvature tensor are bounded.
Under additional mild assumptions, we also know that the injectivity radius is
bounded from below by a positive constant, \ie\ $(M_{0},g_{0})$ is of bounded
geometry. (A {\em manifold with bounded geometry} is a Riemannian manifold
with positive injectivity radius and with bounded covariant derivatives of the
curvature tensor, see \cite{ShubinBook} and references therein). A useful
property is that all geometric operators on $M_0$ that are associated to a
metric on $A$ are $\VV$-differential operator (\ie in $\DiffV{m}(M)$
\cite{aln1}).

On a Riemannian manifold $M_0$ with a Lie structure at infinity $(M, \VV)$,
$\VV = \Gamma(M, A)$, the exponential map $\exp_p : T_pM_0 \to M_0$ is
well-defined for all $p \in M_0$ and extends to a differentiable map $\exp_p :
A_p \to M$ depending smoothly on $p \in M$. A convenient way to introduce the
exponential map is via the geodesic spray, as done in \cite{aln1}. A related
phenomenon is that any vector field $X \in \Gamma(A)$ is integrable, which is
a consequence of the compactness of $M$.  The resulting diffeomorphism of
$M_0$ will be denoted $\psi_X$.

\begin{proposition}\label{prop.bdry.depth}\
  Let $F_0$ be an open face of $M$ and $X\in\Gamma(M; A)$.  Then the
  diffeomorphism $\psi_X$ maps $F_0$ to itself.
\end{proposition}

\begin{proof}
  This follows right away from the assumption that all vector fields in $\VV$
  are tangent to all faces \cite{aln1}.
\end{proof}

\section{Kohn-Nirenberg quantization and pseudodifferential
operators\label{sec.KN}}

{\em Throughout this section $M_0$ will be a fixed manifold with Lie structure
  at infinity $(M, \VV)$ and $\VV := \Gamma(A)$. We shall also fix a metric
  $g$ on $A \to M$, which induces a metric $g_0$ on $M_0$.}  We are going to
introduce a pseudodifferential calculus on $M_0$ that microlocalizes the
algebra of $\VV$-differential operators $\DiffV{*}(M_0)$ on $M$ given by the
Lie structure at infinity.

\subsection{Riemann-Weyl fibration}
Fix a Riemannian metric $g$ on the bundle $A$, and let $g_{0}=g|_{M_{0}}$ be
its restriction to the interior $M_{0}$ of $M$. We shall use this metric to
trivialize all density bundles on $M$. Denote by $\pi : TM_0 \to M_0$ the
natural projection.  Define
\begin{equation}\label{eq.Phi}
  \Phi:TM_{0}\longrightarrow M_{0}\times M_{0}, \quad
  \Phi(v) := (x, \exp_{x}(-v)), \; x = \pi(v).
\end{equation}
Recall that for $v \in T_{x}M$ we have $\exp_{x}(v)=\gamma_{v}(1)$ where
$\gamma_{v}$ is the unique geodesic with $\gamma_{v}(0)=\pi(v)=x$ and
$\gamma_{v}'(0)=v$.  It is known that there is an open neighborhood $U$ of the
zero-section $M_{0}$ in $TM_{0}$ such that $\Phi|_{U}$ is a diffeomorphism
onto an open neighborhood $V$ of the diagonal $M_0=\Delta_{M_{0}}\subseteq
M_{0}\times M_{0}$.

To fix notation, let $E$ be a real vector space together with a metric or a
vector bundle with a metric. We shall denote by $(E)_r$ the set of all vectors
$v$ of $E$ with $|v|< r$.

{\em We shall also assume from now on that $r_0$, the injectivity radius of
  $(M_{0},g_{0})$, is positive.} We know that this is true under some
additional mild assumptions and we conjectured that the injectivity radius is
always positive \cite{aln1}.  Thus, for each $0 < r \le r_{0}$, the
restriction $\Phi|_{(TM_{0})_{r}}$ is a diffeomorphism onto an open
neighborhood $V_{r}$ of the diagonal $\Delta_{M_{0}}$. It is for this reason
that we need the positive injectivity radius assumption.

We continue, by slight abuse of notation, to write $\Phi$ for that
restriction.  Following Melrose, we shall call $\Phi$ the {\em Riemann-Weyl
  fibration}.  The inverse of $\Phi$ is given by $$
M_{0}\times M_{0}\supseteq
V_{r}\ni(x,y)\longmapsto (x,\tau(x,y))\in (TM_{0})_{r}\,, $$
where
$-\tau(x,y)\in T_{x}M_{0}$ is the tangent vector at $x$ to the shortest
geodesic $\gamma : [0,1] \to M$ such that $\gamma(0) = x$ and $\gamma(1) = y$.

\subsection{Symbols and conormal distributions}
Let $\pi: E\rightarrow M$ be a smooth vector bundle with orthogonal metric
$g$.  Let
\begin{equation}\label{eq.<>}
    \<\xi\> := \sqrt{1 + g(\xi,\xi)}.
\end{equation}
We shall denote by $S^m_{1,0}(E)$ the symbols of type $(1,0)$ in H\"ormander's
sense \cite{hor3}.  Recall that they are defined, in local coordinates, by the
standard estimates $$
|\pa_x^\alpha \pa_\xi^\beta a(\xi)| \le
C_{K,\alpha,\beta} \<\xi\>^{m - |\beta|},\quad \pi(\xi) \in K, $$
where $K$ is
a compact trivializing subset (\ie\ $\pi^{-1}(K) \simeq K \times \RR^n$) and
$\alpha$ and $\beta$ are multi-indices. If $a \in S^m_{1, 0}(E)$, then its
image in $S^m_{1, 0}(E)/S^{m-1}_{1, 0}(E)$ is called the {\em principal
  symbol} of $a$ and denoted $\sia{m}(a)$. A symbol $a$ will be called {\em
  homogeneous of degree $\mu$ if} $a(x, \lambda \xi) = \lambda^\mu a(x, \xi)$
for $\lambda > 0$ and $|\xi|$ and $|\lambda \xi|$ large.  A symbol $a \in
S^m_{1, 0}(E)$ will be called {\em classical} if there exist symbols $a_k \in
S^{m-k}_{1, 0}(E)$, homogeneous of degree $m - k$, such that $a -
\sum_{j=0}^{N-1} a_k \in S^{m-N}_{1,0}(E)$.  Then we identify $\sia{m}(a)$
with $a_0$.  (See any book on pseudodifferential operators or the
corresponding discussion in \cite{alnv}.)

We now specialize to the case $E = A^*$, where $A \to M$ is the vector bundle
such that $\VV =\Gamma(M, A)$.  Recall that we have fixed a metric $g$ on $A$.
Let $\pi : A \rightarrow M$ and $\overline{\pi} : A^* \rightarrow M$ be the
canonical projections.  Then the inverse of the Fourier transform $\cF_{{\rm
    fiber}}^{-1}$, along the fibers of $A^{*}$ gives a map
\begin{equation}\label{eq.Finv}
  \cF_{{\rm fiber}}^{-1}: S^{m}_{1,0}(A^{*})
  \longrightarrow\cmun(A) := \cunc(A)' \,, \quad \<\cF_{{\rm
      fiber}}^{-1}a, \phi \> := \<a, \cF_{{\rm fiber}}^{-1}\phi \>,
\end{equation}
where $a \in S^{m}_{1,0}(A^{*})$, $\phi$ is a smooth, compactly supported
function, and
\begin{equation}\label{eq.FI}
    \cF_{{\rm fiber}}^{-1} (\phi)(\xi) := (2\pi)^{-n}
    \int_{\overline{\pi}(\zeta) = \pi(\xi)} e^{i \<\xi, \zeta\>}
    \phi(\zeta) \,d\zeta.
\end{equation}
Then $I^{m}(A, M)$ is defined as the image of $S^m_{1,0}(A^*)$ through the
above map. We shall call this space the space of distributions on $A$ {\em
  conormal} to $M$. The spaces $I^{m}(TM_{0},M_{0})$ and $I^{m}(M_{0}^{2},
\Delta_{M_{0}}) = I^{m}(M_{0}^{2}, M_{0})$ are defined similarly. In fact,
these definitions are special cases of the following more general definition.
Let $X \subset Y$ be an embedded submanifold of a manifold with corners~$Y$.
On a small neighborhood $V$ of $X$ in $Y$ we define a structure of a vector
bundle over $X$, such that $X$ is the zero section of~$V$. As a bundle $V$ is
isomorphic to the normal bundle of $X$ in $Y$.  Then we define the space of
{\em distributions on $Y$ that are conormal of order $m$ to $X$,} denoted
$I^m(Y, X)$, to be the space of distributions on $M$ that are smooth on $Y
\smallsetminus X$ and, that are, in a tubular neighborhood $V \to X$ of $X$ in
$Y$, the inverse Fourier transforms of elements in $S^m(V^*)$ along the fibers
of $V \to X$.  For simplicity, we have ignored the density factor.  For more
details on conormal distributions we refer to \cite{fio, hor3, sim} and the
forthcoming book \cite{MelroseCorners} (for manifolds with corners).

The main use of spaces of conormal distributions is in relation to
pseudodifferential operators. For example, since we have
\begin{equation*}
    I^{m}(M_{0}^{2}, M_{0}) \subseteq
    \cmun(M_{0}^{2}) := \cunc(M_0^2)'\,,
\end{equation*}
we can associate to a distribution in $K \in I^{m}(M_{0}^{2},M_{0})$ a
continuous linear map $T_{K} : \cunc(M_{0})\rightarrow\cmun(M_{0}) :=
\cunc(M_{0})',$ by the Schwartz kernel theorem. Then a well known result of
H\"ormander \cite{fio, hor3} states that $T_K$ is a pseudodifferential
operator on $M_0$ and that all pseudodifferential operators on $M_0$ are
obtained in this way, for various values of $m$. This defines a map
\begin{equation}\label{eq.def.TK}
    T : I^m(M_0^2, M_0)  \to  \Hom(\cunc(M_{0}), \cmun(M_{0})).
\end{equation}

Recall now that $(A)_r$ denotes the set of vectors of norm $< r$ of the vector
bundle $A$.  We agree to write $I^{m}_{(r)}(A, M)$ for all $k\in I^{m}(A, M)$
with $\supp\, k\subseteq (A)_{r}$. The space $I^{m}_{(r)}(TM_{0}, M_{0})$ is
defined in an analogous way. Then restriction defines a map
\begin{equation}\label{eq.def.R}
    \cR: I^{m}_{(r)}(A, M) \longrightarrow I^{m}_{(r)}(TM_{0}, M_{0}).
\end{equation}

Recall that $r_0$ denotes the injectivity radius of $M_0$ and that we assume
$r_0 > 0$. Similarly, the Riemann--Weyl fibration $\Phi$ of Equation
\eqref{eq.Phi} defines, for any $0 < r \le r_0$, a map
\begin{equation}\label{eq.def.Phi_*}
    \Phi_*: I^{m}_{(r)}(TM_0, M_0) \to I^{m}(M_0^2, M_0).
\end{equation}

We shall also need various subspaces of conormal distributions, 
which we shall denote by including a subscript as follows:
\begin{itemize}
\item ``cl'' to designate the distributions that are ``classical,'' in
the sense that they correspond to classical pseudodifferential
operators.
\item ``c'' to denote distributions that have compact support.
\item ``pr'' to indicate operators that are properly supported or
distributions that give rise to such operators.
\end{itemize}
For instance, $I_c^m(Y, X)$ denote
the space of compactly supported conormal distributions, so
that $I^{m}_{(r)}(A, M) = I^{m}_c((A)_{r}, M)$.
Occasionally, we shall use the double subscripts ``cl,pr'' and
``cl,c.'' Note that ``c'' implies ``pr''.  .

\subsection{Kohn-Nirenberg quantization\label{subsec.KN}}\ For notational
simplicity, we shall use the metric $g_{0}$ on $M_0$ (obtained from the metric
on $A$) to trivialize the half-density bundle $\oeh(M_{0})$. In particular, we
identify $\cunc(M_{0}, \oeh)$ with $\cunc(M_{0})$. Let $0 < r \le r_{0}$ be
arbitrary. Each smooth function $\chi$, with $\chi = 1$ close to $M \subseteq
A$ and support contained in the set $(A)_{r}$, induces a map $q_{\Phi,\chi} :
S^{m}_{1,0} (A^{*})\longrightarrow I^m(M_0^2, M_0)$,
\begin{equation}\label{defquant}
    q_{\Phi, \chi}(a) := \Phi_{*}\left(\cR \left(\chi \cF_{{\rm
    fiber}}^{-1}(a) \right) \right).
\end{equation}
Let $a_{\chi}(D)$ be the operator on $M_0$ with distribution kernel
$q_{\Phi,\chi}(a)$, defined using the Schwartz kernel theorem, \ie
$a_{\chi}(D) : = T \circ q_{\Phi, \chi}(a)$ .  Following Melrose, we call the
map $q_{\Phi,\chi}$ the {\em Kohn-Nirenberg quantization map}. It will play an
important role in what follows.

For further reference, let us make the formula for the induced operator
$a_\chi(D) : \cunc(M_{0})\rightarrow\cunc(M_{0})$ more explicit. Neglecting
the density factors in the formula, we obtain for $u\in\cunc(M_{0})$
\begin{equation}\label{eq.12}
  a_\chi(D)u(x) = \int_{M_{0}} (2\pi )^{-n}
  \int_{T^{*}_{x}M_{0}} e^{i\tau(x,y)\cdot\eta}
  \chi(x,\tau(x,y)) a(x,\eta)u(y)\, d \eta \,dy \,.
\end{equation}
Specializing to the case of Euclidean space $M_{0} = \RR^{n}$ with the
standard metric we have $\tau(x,y) = x - y$, and hence
\begin{equation}\label{weylq}
  a_\chi(D)u(x) = (2\pi )^{-n} \int_{\RR^{n}}\int_{\RR^{n}}
  e^{i(x-y)\eta} \chi(x,x-y) a(x,\eta)u(y) \,d \eta \,dy\,,
\end{equation}
\ie\ the well-known formula for the Kohn-Nirenberg-quantization on $\RR^{n}$,
if $\chi = 1$. The following lemma states that, up to regularizing operators,
the above quantization formulas do not depend on $\chi$.

\begin{lemma}\label{lemma.cont.R}\
  Let $0 < r \le r_0$.  If $\chi_1$ and $\chi_2$ are smooth functions with
  support $(A)_{r}$ and $\chi_j = 1$ in a neighborhood of $M\subseteq A$, then
  $(\chi_1 - \chi_2)\cF_{{\rm fiber}}^{-1}(a)$ is a smooth function, and hence
  $a_{\chi_1}(D) - a_{\chi_2}(D)$ has a smooth Schwartz kernel. Moreover, the
  map $S^{m}_{1,0}(A^*) \rightarrow \CI(A)$ that maps $a\in S^m_{1,0}(A^*)$ to
  $(\chi_1 - \chi_2)\cF_{{\rm fiber}}^{-1}(a)$ is continuous, where the right
  hand side is endowed with the topology of uniform $C^\infty$-convergence on
  compact subsets.
\end{lemma}

\begin{proof}\ Since the singular supports of $\chi_1\cF_{{\rm
      fiber}}^{-1}(a)$ and $\chi_2\cF_{{\rm fiber}}^{-1}(a)$ are contained in
  the diagonal $\Delta_{M_{0}}$ and $\chi_{1}-\chi_{2}$ vanishes there, we
  have that $(\chi_1 - \chi_2)\cF_{{\rm fiber}}^{-1}(a)$ is a smooth function.

  To prove the continuity of the map $S^{m}_{1,0}(A^*) \ni a \mapsto (\chi_1 -
  \chi_2)\cF_{{\rm fiber}}^{-1}(a) \in \CI(A)$, it is enough, using a
  partition of unity, to assume that $A \to M$ is a trivial bundle. Then our
  result follows from the standard estimates for oscillatory integrals (\ie by
  formally writing $|v|^2\int e^{i \<v, \xi\>} a(\xi) d\xi = -\int (\Delta_\xi
  e^{i \<v, \xi\>}) a(\xi) d\xi$ and then integrating by parts, see
  \cite{hor3, NagelStein, Taylor1, Taylor2} for example).
\end{proof}

We now verify that the quantization map $q_{\Phi, \chi}$, Equation
\eqref{defquant}, gives rise to pseudodifferential operators.

\begin{lemma}\label{basqu}\
  Let $r \le r_{0}$ be arbitrary. For each $a\in S_{1,0}^{m}(A^{*})$ and each
  $\chi\in\cunc((A)_{r})$ with $\chi=1$ close to $M\subseteq A$, the
  distribution $q_{\Phi,\chi}(a)$ is the Schwartz-kernel of a
  pseudodifferential operator $a_\chi(D)$ on $M_{0}$, which is properly
  supported if $r < \infty$ and has principal symbol $\sia{\mu}(a) \in S^m_{1,
    0}(E)/S^{m-1}_{1, 0}(E)$. If $a \in S_{cl}^{\mu}(A^*)$, then $a_\chi(D)$
  is a classical pseudodifferential operator.
\end{lemma}

\begin{proof}\ Denote also by $\chi : I^{m}(TM_0, M_0) \to I^{m}_{(r)}(TM_0,
  M_0)$ the ``multiplication by $\chi$'' map. Then
\begin{equation}
  a_\chi(D) = T \circ \Phi_* \circ \cR \circ \chi \circ \cF_{{\rm
      fiber}}^{-1}(a) := T_{ \Phi_*(\cR( \chi \cF_{{\rm fiber}}^{-1}(a)
    ))} = T \circ q_{\Phi, \chi} (a)
\end{equation}
where $T$ is defined in Equation \eqref{eq.def.TK}. Hence $a_\chi(D)$ is a
pseudodifferential operator by the H\"ormander's result mentioned above
\cite{fio, hor3} (stating that the distribution conormal to the diagonal are
exactly the Schwartz kernels of pseudodifferential operators.  Since $\chi
\cR(a)$ is properly supported, so will be the operator $a_\chi(D)$).

For the statement about the principal symbol, we use the principal symbol map
for conormal distributions \cite{fio, hor3}, and the fact that the restriction
of the anchor $A\rightarrow TM$ to the interior $A|_{M_{0}}$ is the identity.
(This also follows from the Equation \eqref{weylq} below.)  This proves our
lemma.
\end{proof}

Let us denote by $\Psi^m(M_0)$ the space of pseudodifferential operators of
order $\le m$ on $M_0$ (no support condition). We then have the following
simple corollary.

\begin{corollary}\label{cor.comp.symb}\
The map $\sigma_{tot} : S_{1,0}^{m}(A^{*}) \to
\Psi^{m}(M_{0})/\Psi^{-\infty}(M_{0})$,
$$
    \sigma_{tot}(a) := a_\chi(D)+\Psi^{-\infty}(M_{0})
$$ is independent of the choice of the function $\chi\in
\CIc((A)_{r})$ used to define $a_\chi(D)$ in Lemma {\rm
\ref{basqu}}.
\end{corollary}

\begin{proof}\ 
This follows right away from the previous lemma, Lemma \ref{basqu}.
\end{proof}

Let us remark that our pseudodifferential calculus depends on more than just
the metric.

\begin{remark}
Non-isomorphic Lie structures at infinity can lead to the same metric
on $M_0$. An example is provided by $\RR^n$ with the standard metric,
which can be obtained either from the radial compactification of
$\RR^n$ with the scattering calculus, or from $[-1, 1]^n$ with the
$b$-calculus.  See Examples \ref{ex1} and the paragraph following it.
The pseudodifferential calculi obtained from these Lie algebra
structures at infinity will be, however, different.
\end{remark}

The above remark readily shows that not all pseudodifferential
operators in $\Psi^{m}(M_{0})$ are of the form $a_\chi(D)$ for some
symbol $a\in S_{1,0}^{m}(A^{*})$, not even if we assume that they are
properly supported, because they do not have the correct behavior at
infinity. Moreover, the space $T \circ q_{\Phi, \chi}
(S_{1,0}^{\infty} (A^{*}) )$ of all pseudodifferential operators of
the form $a_\chi(D)$ with $a\in S_{1,0}^{\infty}(A^{*})$ is not closed
under composition. In order to obtain a suitable space of
pseudodifferential operators that is closed under composition, we are
going to include more (but not all) operators of order $-\infty$ in
our calculus.

Recall that we have fixed a manifold $M_0$, a Lie structure at
infinity $(M,A)$ on $M_0$, and a metric $g$ on $A$ with injectivity 
radius $r_{0}>0$. Also, recall that any $X \in \Gamma(A) \subset \VV_b$ 
generates a global flow $\Psi_{X} : \RR\times M\rightarrow M$. 
Evaluation at $t=1$ yields a diffeomorphism $\Psi_{X}(1,\cdot) : M \rightarrow M$,
whose action on functions is denoted 
\begin{equation}\label{eq.def.psiX}
    \psi_X : \CI(M) \to \CI(M).
\end{equation}

We continue to assume that the injectivity radius $r_0$ of our fixed
manifold with a Lie structure at infinity $(M, \VV)$ is strictly
positive.

\begin{definition}\label{def.psi.MA}\ Fix $0< r < r_0$ and
$\chi \in \cunc((A)_{r})$ such that $\chi = 1$ in a neighborhood of $M
\subseteq A$.  For $m \in \RR$, the space $\Psi_{1,0,\VV}^{m}(M_{0})$
of {\em pseudodifferential operators generated by the Lie structure at
infinity $(M,A)$} is the linear space of operators $\cunc(M_0)
\rightarrow \cunc(M_{0})$ generated by $a_\chi(D)$, $a \in
S_{1,0}^m(A^*)$, and $b_{\chi}(D)\psi_{X_1}\ldots \psi_{X_k}$, $b \in
S^{-\infty}(A^*)$ and $X_j \in \Gamma(A)$, $\forall j$.

Similarly, the space $\Psi_{cl,\VV}^{m}(M_{0})$ of {\em classical
pseudodifferential operators generated by the Lie structure at
infinity $(M,A)$} is obtained by using classical symbols $a$ in the
construction above.
\end{definition}

It is implicit in the above definition that the spaces
$\Psi_{1,0,\VV}^{-\infty}(M_{0})$ and $\Psi_{cl,\VV}^{-\infty}(M_{0})$
are the same. They will typically be denoted by
$\Psi^{-\infty}_{\VV}(M_0)$. As usual, we shall denote
\begin{equation*}
    \Psi_{1,0,\VV}^{\infty}(M_{0}) := \cup_{m \in \ZZ}
    \Psi_{1,0,\VV}^{m}(M_{0}) \quad \text{and} \quad
    \Psi_{cl,\VV}^{\infty}(M_{0}) := \cup_{m \in \ZZ}
    \Psi_{cl,\VV}^{m}(M_{0}).
\end{equation*}

At first sight, the above definition depends on the choice of the
metric $g$ on $A$. However, we shall soon prove that this is not the
case.

As for the usual algebras of pseudodifferential operators, we have the
following basic property of the principal symbol.

\begin{proposition}\label{prop.princ.symb}\
The principal symbol establishes isomorphisms
\begin{equation}\label{eq.def.princ}
    \sigma^{(m)} : \Psi_{1,0,\VV}^{m}(M_{0}) /
    \Psi_{1,0,\VV}^{m-1}(M_{0}) \to S^m_{1,0}(A^*) /
    S^{m-1}_{1,0}(A^*)
\end{equation}
and
\begin{equation}\label{eq.def.princ'}
    \sigma^{(m)} : \Psi_{cl,\VV}^{m}(M_{0}) /
    \Psi_{cl,\VV}^{m-1}(M_{0}) \to S^m_{cl}(A^*) / S^{m-1}_{cl}(A^*).
\end{equation}
\end{proposition}

\begin{proof}\
This follows from the classical case of the spaces $\Psi^m(M_0)$ using
also Lemma \ref{basqu}.
\end{proof}

\section{The product\label{sec.product}}

We continue to denote by $(M, \VV)$, $\VV = \Gamma(A)$, a fixed
manifold with a Lie structure at infinity and with positive
injectivity radius.  In this section we want to show that the space
$\Psi_{1,0,\VV}^{\infty}(M_{0})$ is an algebra (\ie it is closed under
multiplication) by showing that it is the homomorphic image of the algebra 
$\Psi_{1, 0}^\infty\GR)$ of pseudodifferential operators on any $d$-connected 
groupoid $\GR$ integrating $A$ (Theorem~\ref{theorem.gr}).

First we need to fix the terminology and to recall some definitions and
constructions involving groupoids.

\subsection{Groupoids}
Here is first an abstract definition that will be made more clear below.
Recall that a small category is a category whose morphisms form
a set. A {\em groupoid} is a small category all of whose morphisms are
invertible. Let $\GR$ denote the set of morphisms and $M$ denote the
set of objects of a given groupoid. Then each $g \in \GR$ will have a
domain $d(g) \in M$ and a range $r(g) \in M$ such that the product
$g_1g_2$ is defined precisely when $d(g_1) = r(g_2)$. Moreover, it
follows that the multiplication (or composition) is associative and
every element in $\GR$ has an inverse. We shall identify the set of
objects $M$ with their identity morphisms via a map $\iota : M \to \GR$. 
One can think then of a groupoid as being a group, except that the
multiplication is only partially defined. By abuse of notation, we
shall use the same notation for the groupoid and its set of morphisms
($\GR$ in this case). An intuitive way of thinking of a groupoid with
morphisms $\GR$ and objects $M$ is to think of the elements of $\GR$
as being arrows between the points of $M$. The points of $M$ will be called 
{\em units}, by identifying an object with its identity morphism. There will be
structural maps $d, r : \GR \to M$, {\em domain and range}, 
$\mu : \{(g, h), d(g)=r(h)\} \to \GR$, {\em multiplication}, 
$\GR \ni g \to g^{-1} \in \GR$, {\em inverse}, and 
$\iota : M \to \GR$ satisfying the usual identities satisfied by
the composition of functions.

A {\em Lie groupoid} is a groupoid $\GR$ such that the
space of arrows $\GR$ and the space of units $M$ are manifolds with
corners, all its structural maps (\ie multiplication, inverse, domain,
range, $\iota$) are differentiable, the domain and range maps (\ie $d$ and 
$r$) are submersions. By the definition of a submersion of manifolds with
corners, the submanifolds $\GR_x := d^{-1}(x)$ and $\GR^x := r^{-1}(x)$ are
smooth (so they have no corners or boundary), for any $x \in M$. Also, 
it follows that that $M$ is an embedded
submanifold of $\GR$.

The {\em $d$--vertical tangent space to $\GR$}, denoted $T_{vert}\GR$,
is the union of the tangent spaces to the fibers of $d : \GR \to M$,
that is
\begin{equation}\label{eq.def.Tvert}
    T_{vert}\GR := \cup_{x \in M} T\GR_x = \ker d_*,
\end{equation}
the union being a disjoint union, with topology induced from the
inclusion $T_{vert} \GR \subset T\GR$. The {\em Lie algebroid of $\GR$},
denoted $A(\GR)$ is defined to be the restriction of the $d$--vertical
tangent space to the set of units $M$, that is, $A(\GR) = \cup_{x \in M}
T_x\GR_x$, a vector bundle over $M$. The space of sections of
$A(\GR)$ identifies canonically with the space of sections of the
$d$-vertical tangent bundle (= $d$-vertical vector fields) that are
{\em right invariant} with respect to the action of $\GR$. It also
implies a canonical isomorphism between the vertical tangent bundle
and the pull-back of $A(\GR)$ via the range map $r : \GR \to M$:
\begin{equation}\label{eq.isom.r}
    r^*A(\GR) \simeq T_{vert}\GR.
\end{equation}
The structure of Lie algebroid on $A(\GR)$ is induced by the Lie
brackets on the spaces $\Gamma(T \GR_x)$, $\GR_x := d^{-1}(x)$.  This
is possible since the Lie bracket of two right invariant vector fields
is again right invariant. The anchor map in this case is given by the
differential of $r$, $r_* : A(\GR) \to TM$.

Let $\GR$ be a Lie groupoid with units $M$, then there is
associated to it a pseudodifferential calculus (or algebra of
pseudodifferential operators) $\Psi^{\infty}_{1,0}(\GR)$, whose
operators of order $m$ form a linear space denoted $\Psi^{m} _{1,0}
(\GR)$, $m \in \RR$, such that $\Psi^m_{1,0} (\GR) \Psi^{m'}_{1,0}
(\GR) \subset \Psi^{m + m'}_{1,0} (\GR)$.  This calculus is defined as
follows:\ $\Psi^{m}_{1,0}(\GR)$ consists of smooth families of
pseudodifferential operators $(P_x)$, $P_x \in \Psi_{1, 0}^m(\GR_x)$, 
$x \in M$, that are right 
invariant with respect to multiplication by elements of $\GR$ and are
``uniformly supported.'' To define what uniformly supported means, let
us observe that the right invariance of the operators $P_x$ implies
that their distribution kernels $K_{P_x}$ descend to a distribution
$k_P \in I^m(\GR, M)$. Then the family $P = (P_x)$ is called {\em
uniformly supported} if, by definition, $k_P$ has compact support. If
$P$ is uniformly supported, then each $P_x$ is properly supported.
The right invariance condition means, for $P = (P_x) \in
\Psi^\infty_{1,0}(\GR)$, that right multiplication $\GR_x \ni g'
\mapsto g'g \in \GR_y$ maps $P_y$ to $P_x$, whenever $d(g) = y$ and
$r(g) = x$. By definition, the evaluation map
\begin{equation}\label{eq.ex}
    \Psi^\infty_{1,0}(\GR) \ni P=(P_x) \mapsto e_z(P) := P_z \in
    \Psi^{\infty}_{1,0}(\GR_z)
\end{equation}
is an algebra morphism for any $z \in M$.
If we require that the operators $P_x$ be classical of order $\mu \in
\CC$, we obtain spaces $\Psi_{cl}^{\mu}(\GR)$ having similar
properties. These spaces were considered in \cite{nwx}.

All results and constructions above remain true for classical
pseudodifferential operators. This gives the algebra
$\Psi_{cl}^\infty(\GR)$ consisting of families $P = (P_x)$ of {\em
classical} pseudodifferential operators satisfying all the previous
conditions.

Assume that the interior $M_{0}$ of $M$ is an invariant subset.
Recall that the so called {\em vector representation} $\piM : \Psi_{1,
0}^{\infty}(\GR) \to \End(\CIc(M_0))$ associates to a
pseudodifferential operator $P$ on $\GR$ a pseudodifferential operator
$\piM(P) : \cunc(M_{0})\rightarrow \cunc(M_{0})$ \cite{LN1}. This
representation $\piM$ is defined as follows.  If $\phi \in \CIc(M_0)$,
$\phi \circ r$ is a smooth function on $\GR$, and we can let the
family $(P_x)$ act along each $\GR_x$ to obtain the function $P(\phi
\circ r)$ on $\GR$ defined by $P(\phi \circ r) \vert_{\GR_{x}} = P_x
(\phi \circ r \vert_{\GR_{x}})$. The fact that $P_x$ is a smooth
family guarantees that $P(\phi \circ r)$ is also smooth. Using then
the fact that $r$ is a submersion, so locally it is a product map, we
obtain that $P(\phi \circ r) = \phi_0 \circ r,$ for some function
$\phi_0 \in \CIc(M_0)$. We shall then let
\begin{equation}\label{eq.def.piM}
    \piM(P) \phi = \phi_0.
\end{equation}
The fact that $P$ is uniformly supported guarantees that $\phi_0$
will also have compact support in $M_0$. A more explicit 
description of $\piM$ in the case of Lie manifolds 
will be obtained in the proof of Theorem \ref{theorem.gr}, more
precisely, Equation \eqref{eq.25}.

A Lie groupoid $\GR$ with units $M$ is said to 
{\em integrate} $A$ if $A(\GR) \simeq A$ as vector bundles over $M$. 
Recall that the groupoid $\GR$ is called {\em $d$--connected} if $\GR_x
:= d^{-1}(x)$ is a connected set, for any $x \in M$. If there exists a
Lie groupoid $\GR$ integrating $A$, then there exists also a $d$--connected 
Lie groupoid with this property. (Just take for each $x$ 
the connected component of $x$ in $\GR_x$.)

Our plan to show that $\Psi_{1, 0, \VV}^\infty(M_0)$ is an algebra,
is then to prove that it is the image under $\piM$ of $\Psi_{1,0}^\infty(\GR)$,
for a Lie groupoid $\GR$ integrating $A$, $\Gamma(M, A) = \VV$. In fact, 
any $d$-connected Lie groupoid will satisfy this, by Theorem 
\ref{theorem.gr}. This
requires the following deep result due to Crainic and 
Fernandes \cite{CrainicFernandes} stating that the Lie algebroids associated
to Lie manifolds are integrable.

\begin{theorem}\label{exists}[Cranic--Fernandes]
Any Lie algebroid arising from a Lie structure at infinity is actually
the Lie algebroid of a Lie groupoid (\ie it is {\em
integrable}).
\end{theorem}

This Theorem should be thought of as an analog of Lie's third
theorem stating  that every
finite dimensional Lie algebra is the Lie algebra of a Lie group.
However, the analog of Lie's theorem for Lie algebroids does not hold:
there are Lie algebroid which are not Lie algebroids to a Lie groupoid
\cite{Mackenzie}.

A somewhat weaker form of the above theorem, 
which is however enough for the proof of Melrose's conjecture was 
obtained in \cite{NistorINT}.

We are now ready to state and prove the main result of this
section. We refer to \cite{LN1} or \cite{nwx} for the concepts and results 
on groupoids and algebras of pseudodifferential operators on groupoids
not explained below or before the statement of this theorem.

\begin{theorem}\label{theorem.gr}\
Let $M_0$ be a manifold with a Lie structure at infinity, 
$(M,\VV)$, $A=A_\VV$, as
above. Also, let $\GR$ be a $d$-connected groupoid with units $M$ and
with $A(\GR) \simeq A$. Then $\Psi_{1,0,\VV}^{m}(M_{0}) = \piM(
\Psi_{1,0}^m(\GR))$ and $\Psi_{cl,\VV}^{m}(M_{0}) = \piM(
\Psi_{cl}^m(\GR))$.
\end{theorem}

\begin{proof}\
We shall consider only the first equality. The case of classical
operators can be treated in exactly the same way.

Here is first, briefly, the idea of the proof. Let $P = (P_x) \in
\Psi^m_{1,0}(\GR)$.  Then the Schwartz kernels of the operators $P_x$ 
form a smooth family of conormal distributions in $I^m(\GR_x^2, \GR_x)$
that descends, by right invariance, to a distribution $k_P \in I_{c}^m(\GR,
M)$ (\ie to a compactly supported distribution on $\GR$, conormal to
$M$) called the {\em convolution kernel of $P$}. The map $P \mapsto k_P$
is an isomorphism \cite{nwx} with inverse 
\begin{equation}\label{eq.T}
  T: I_c^m(\GR, M) \to
  \Psi^m_{1,0}(\GR).
\end{equation}
Fix a metric on $A \to M$. The resulting exponential map (reviewed
below) then gives rise, for $r > 0$ small enough, to an open embedding
\begin{equation}\label{eq.def.alpha}
    \alpha : (A)_r \to \GR,
\end{equation}
which is a diffeomorphism onto its image. This diffeomorphism then
gives rise to an embedding
\begin{equation}\label{eq.emb.Ar}
    \alpha_* : I^m_{(r)}(A, M) := I^m_c((A)_r, M) \to I^m(\GR,M)
\end{equation}
such that for each $\chi$ as above
\begin{equation}\label{eq.eq.Ar}
    \piM\big(\alpha_*( \chi \cF_{{\rm fiber}}^{-1}(a) )\big) =
    a_{\chi}(D) \in \Psi^m(M_0).
\end{equation}
This will allow us to show that $\piM(\Psi_{1, 0}^m(\GR))$ contains the
linear span of all operators $P$ of the form $P = a_\chi(D)$, $a \in
S^m_{1,0}(A^*)$, $m \in \ZZ$ fixed. This reduces the problem to
verifying that
\begin{equation}\label{eq.prof.reg}
    \piM(\Psi^{-\infty}(\GR)) = \Psi_{\VV}^{-\infty}(M_0).
\end{equation}
Using a partition of unity, this in turn will be reduced to Equation
\eqref{eq.eq.Ar}. Now let us provide the complete details.

Let $\GR_x^x := d^{-1}(x) \cap r^{-1}(x)$, which is a group for any $x
\in M_0$, by the axioms of a groupoid. Then $\GR_x^x \simeq \GR_y^y$
whenever there exists $g \in \GR$ with $d(g) = x$ and $r(g) = y$
(conjugate by $g$).  We can assume, without loss of generality, that $M$ 
is connected. Let $\Gamma := \GR_x^x$, for some fixed $x \in M_0$.
Our above informal description of the proof can be
conveniently formalized and visualized using the following diagram
whose morphisms are defined below:

\begin{equation*}
\label{comm.diagram}
\xymatrix{
S^{m}_{1,0}(A)      \ar[r]^{\mF^{-1}_{\rm fiber}\;}             &
I^{m}(A, M)         \ar[r]^{\chi}               &
I^{m}_{(r)}(A, M)   \ar[d]_{\alpha_*} \ar[r]^{\mR\;}    &
I^{m}_{(r)}(TM_0, M_0)  \ar[d]^{\Phi_*}
\\
                                &
\Psi^{m}_{1,0}(\GR)         \ar[r]^{\cong\,} \ar[d]_{e_x}           &
I_{c}^{m}(\GR, M)   \ar[r]^{l_{*}} \ar[d]_{\mu_{1}^{*}}     &
I^{m}(M_0^2, M_0)   \ar@{=}[d]              \\
                                &
\Psi^{m}_{pr}(\GR_x)^\Gamma        \ar[r]^{\cong\,}        &
I^{m}_{pr}(\GR_x^2, \GR_x)^\Gamma
\ar[d]_{\cong\,} \ar[r]^{\tilde{r}_*}   &
I^{m}(M_0^2, M_0)   \ar[d]^{\cong\,}                \\
                                &
                                &
\Psi^{m}_{pr}(\GR_x)^\Gamma \ar[r]^{r_*}                &
\Psi^{m}_{1, 0}(M_0)                      }
\end{equation*}
\vspace*{1mm}

We now define the morphisms appearing in the above diagram in such a
way that it will turn out to be a commutative diagram. The bottom
three rectangles will trivially turn out to be commutative. Recall
that the index ``pr'' means ``properly supported.''

First, recall that the maps $\mF^{-1}_{\rm fiber}$ (the fiberwise
inverse Fourier transform), $\chi$ (the multiplication by the cut-off
function $\chi$), $\mR$ (the restriction map), $\Phi_*$ (induced by
the inverse of the exponential map), and $e_x$ (evaluation of the
family $(P_y)$ at $y = x$) have already been defined.

We let $\mu_1(g',g) = g'g^{-1}$ and $\mu_1^*$ be the map induced at
the level of kernels by $\mu_1$ by pull-back (which is seen to be
defined in this case because $\mu_1$ is a submersion and its range is
transverse to $M$).

The four isomorphisms not named are the ``$T$ isomorphisms'' and their
inverses defined in various places earlier (identifying spaces of
conormal distributions with spaces of pseudodifferential
operators). More precisely, the top isomorphism is from \cite{nwx} and
all the other isomorphisms are the canonical identifications between
pseudodifferential operators and distributions on product spaces that
are conormal to the diagonal (via the Schwartz kernels). In fact, the
top isomorphism $T$ is completely determined by the requirement that
the left-most square (containing $e_x$) be commutative.

It is a slightly more difficult task to define $r_*$. We shall have
to make use of a minimum of groupoid theory. Let $y \in M$ be
arbitrary for a moment. 
Since the Lie algebra of
$\GR_y^y$ is isomorphic to the kernel of the anchor map $\varrho :
A(\GR)_y \to T_y M$, we see that $\GR_y^y$ is a discrete group if, and
only if, $y \in M_0$.
Then
\begin{equation*}
        r_* : T_y \GR_y = A(\GR)_y \to T_y M_0
\end{equation*}
is an isomorphism, if and only if, $y \in M_0$. 

Let $x \in M_0$ be our fixed point.  Then $r : \GR_{x} \to M_0$
is a surjective local diffeomorphism. 
Also $\Gamma := \GR_x^x$ acts freely on $\GR_x$ and 
$\GR_x /\Gamma = M_0$. 
Hence $r : \GR_x \to M_0$ is a covering map with group $\Gamma$,
and $\CI(\GR_x)^{\Gamma} = \CI(M_0)$. Let $P = (P_y) \in
\Psi_{1, 0}^m(\GR)$. Since $P_x$
is $\Gamma$-invariant and properly supported, the map $P_x :
\CI(\GR_x) \to \CI(\GR_x)$, descends to a map $\CI(M_0) \to \CI(M_0)$,
which is by definition $r_*(P)$. More precisely, if $\phi$ is a smooth
function on $M_0$, then $\phi \circ r$, is a
$\Gamma$-invariant function on $\GR_x$. Hence $P (\phi \circ r)$ is
defined (because $P$ is properly supported) and is also $\Gamma$-invariant. 
Thus there exists a function $\phi_0 \in \CI(M_0)$ such
that $P(\phi \circ r) = \phi_0 \circ r$. The operator $r_*(P)$ is then
given by $r_*(P)\phi := \phi_0$. This definition of $r_*$ provides us with the following
simpler definition of the vector representation $\piM$:
\begin{equation}\label{eq.25}
    \piM(P) = r_*(e_x(P)).
\end{equation}

We also obtain that
\begin{equation}\label{eq.piM}
    \piM \circ T = r_* \circ e_x \circ T = r_* \circ T \circ \mu_1^*,
\end{equation}
by the commutativity of the left-most rectangle.

The commutativity of the bottom rectangle completely determines the
morphism $\tilde{r}_*$. However, we shall also need an explicit description
of this map. This can be obtained as follows. Recall that $\GR_x$ is a
covering of $M_0$ with group $\Gamma := \GR_x^x$.  Hence we can identify $I^m_{pr}(\GR_x^2, \GR_x)^\Gamma$ with
$I^m_{pr}((\GR_x^2)^\Gamma, \GR_x^\Gamma)$. The map $\tau :
(\GR_x^2)/\Gamma \to M_0^2$ is also a covering map. This allows us to
identify a distribution with small support in $(\GR_x^2)/\Gamma$ with
a distribution with support in a small subset of $M_0^2$. These
identifications then extend by summation along the fibers of $\tau :
(\GR_x^2)/\Gamma \to M_0^2$ to define a distribution $\tau_*(u) \in
\mathcal D'(M_0^2)$, for any distribution $u$ on $(\GR_x^2)/\Gamma$
whose support intersects only finitely
many components of $\tau^{-1}(U)$, for any connected locally trivializing open
set $U \subset M_0$. The morphism $\tilde{r}_*$ identifies then with~$\tau_*$.  
Also, observe for later use that
\begin{equation}\label{eq.lrm}
    \tau(g', g) = (r(g'), r(g)) = (r(g'g^{-1}), d(g'g^{-1})) = ( r(
    \mu_1(g', g) ), d(\mu_1(g', g) ) ).
\end{equation}

Next, we must set $l_* := \tilde{r}_* \circ \mu_1^*$, by the
commutativity requirement. For this morphism we have a similar,
but simpler, description of $l_*(u)$. Namely, $l_*(u)$ is
obtained by first restricting a
distribution $u$ to $d^{-1}(M_0) = r^{-1}(M_0)$ and then by applying
to this restriction the push-forward map defined by $(d, r) : d^{-1}(M_0) \to M_0^2$
(that is, we sum over open sets in $\GR_x$ covering sets in the base $M_0^2$).
Equation \eqref{eq.lrm} guarantees that this alternative description of $l_*$
satisfies $l_* := \tilde{r}_* \circ \mu_1^*$.

To define $\alpha_*$, recall that we have fixed a metric on $A$. This metric then
lifts via $r : \GR \to M$ to $T_{vert}\GR \simeq r^*A(\GR)$, by
Equations \eqref{eq.def.Tvert} and \eqref{eq.isom.r}. The induced
metrics on the fibers of $\GR_y$, $y \in M$, give rise, using the 
(geodesic) exponential map, to maps
\begin{equation*}
    A_y \simeq A(\GR)_y = T_y \GR_y\to \GR_y.
\end{equation*}
These maps give rise to an application $(A)_r \to \GR$, which, by the Inverse
mapping theorem, is seen to be a diffeomorphism onto its image. It,
moreover, sends the zero section of $A$ to the units of $\GR$. Then
$\alpha_*$ is the resulting map at the level of conormal
distributions.  (Note that $\GR_y$ is complete.)

We have now completed the definition of all morphisms in our 
diagram. To prove that our diagram is commutative, it remains to prove
that
\begin{equation*}
    l_* \circ \alpha_* = \Phi_* \circ \mR.
\end{equation*}
This however follows from the above description of the map $l_*$,
since $(d,r)$ is injective on $\alpha((A)_r)$ and $r : \GR_x \to M_0$
is an isometric covering, thus preserving the exponential maps.

The commutativity of the above diagram finally shows that
\begin{multline}\label{eq.surj}
    a_\chi(D) := T \circ q_{\Phi, \chi}(a) = T \circ \Phi_* \circ \mR
    \circ \chi \circ \mF^{-1}_{\rm fiber} (a) \\ = \piM \circ T
    \circ \alpha_* \circ \chi \circ \mF^{-1}_{\rm fiber} (a) =
    \piM(Q),
\end{multline}
where $Q = T \circ \alpha_* \circ \chi \circ \mF^{-1}_{\rm fiber} (a)$
and $a \in S^m_{1,0}(A^*)$. Thus every operator of the form
$a_\chi(D)$ is in the range of $\piM$.

Let us notice for the rest of our argument that the definition of the
vector representation $\piM$ can be extended by the same formula to arbitrary
right invariant families of operators $P = (P_x)$, $P_x : \CI(\GR_x)
\to \CI(\GR_x)$, such that the induced operator $P : \CIc(\GR) \to
\cup \CI(\GR_x)$ has range in $\CIc(\GR)$. We shall use this in the
following case.  Let $X \in \VV$. Then $X$ defines by integration a
diffeomorphism of $M$, see Equation \eqref{eq.def.psiX}. Let $\tilde
X$ be its lift to a $d$-vertical vector field on $\GR$ (\ie on each
$\GR_x$ we obtain a vector field, and this family of vector fields is
right invariant).  A result
from \cite[Appendix]{Spencer} (see also \cite{NistorINT}) then shows
that $\tilde X$ can be integrated to a global flow. Let us denote by
$\tilde{\psi}_{X}$ the {\em family} of diffeomorphisms of each $\GR_x$
obtained in this way, as well as their action on functions. 
It follows then from the definition that
\begin{equation}\label{eq.pi.psi}
    \piM (\tilde{\psi}_X) = \psi_X.
\end{equation}

The Equations \eqref{eq.surj} and \eqref{eq.pi.psi} then give
\begin{equation}\label{eq.vect}
    \piM(Q \tilde{\psi}_{X_1} \ldots \tilde{\psi}_{X_n}) = a_\chi(D)
    \psi_{X_1} \ldots \psi_{X_n} \in \Psi^{-\infty}_{1,0,\VV}(M_0),
\end{equation}
for any $a \in S^{-\infty}(A^*)$ and $Q = T \circ \alpha_* \circ \chi
\circ \mF^{-1}_{\rm fiber} (a)$. Also $Q \tilde{\psi}_{X_1} \ldots
\tilde{\psi}_{X_n} \in \Psi^{-\infty}(\GR)$, since the product of a
regularizing operator with the operator induced by a diffeomorphism is
regularizing. We have thus proved that $\piM(\Psi^{m}_{1,0}(\GR))
\supset \Psi^m_{1,0, \VV}(M_0)$.  Let us now prove the opposite
inclusion, that is that $\piM(\Psi^{m}_{1,0}(\GR)) \subset \Psi^m_{1,0,
\VV}(M_0)$.

Let $Q \in \Psi_{1,0}^m(\GR)$ be arbitrary and let $b = T^{-1}(Q)$.
Let $\chi_0$ be a smooth function on $\GR$ that is equal to 1 in a
neighborhood of $M$ in $\GR$ and with support in $\alpha((A)_r)$ and
such that $\chi = 1$ on the support of $\chi_0 \circ \alpha$. Then
$b_0 := \chi_0 b$ is in the range of $\alpha_* \circ \chi \circ
\mF^{-1}_{\rm fiber}$, because any distribution $u \in I^m_{(r)}(A,M)$
is in the range of $\mF^{-1}_{\rm fiber}$, if $r < \infty$. Then the
difference $b - b_0$ is smooth. Because $\GR$ is $d$-connected, we can
use a similar construction as the one used to define $b_0$ and a partition 
of unity argument to obtain that
\begin{equation}\label{eq.33}
        T(b - b_0) = \sum_{j = 1}^l T(b_j) \tilde{\psi}_{X_{j1}}
        \ldots \tilde{\psi}_{X_{jn}}
\end{equation}
for some distributions $b_j \in \chi I^{-\infty}_{(r)}(A,M)$ and
vector fields $X_{jk} \in \VV$.  Let $a_j$ be such that $b_j =
\alpha_* \circ \chi \circ \mF^{-1}_{\rm fiber} (a_j)$, for $a_0 \in
S_{1,0}^{m}(A^*)$ and $a_j \in S_{1,0}^{-\infty}(A^*)$, if $j >
0$. Then Equations \eqref{eq.surj} and \eqref{eq.vect} show that
\begin{equation}
    \piM(Q) = a_0(D) + \sum_{j = 1}^l a_j(D) \psi_{X_{j1}} \ldots
        \psi_{X_{jn}} \in \Psi^m_{1,0,\VV}(M_0).
\end{equation}
We have thus proved that $\piM(\Psi_{1,0}^m(\GR) ) =
\Psi^m_{1,0,\VV}(M_0),$ as desired. This completes our proof.
\end{proof}

Since the map $\piM$ respects adjoints:\ $\piM(P^*) = \piM(P)^*$,
\cite{LN1}, we obtain the following corollary.

\begin{corollary} The algebras $\Psi_{1, 0, \VV}^{\infty}(M_0)$
and $\Psi_{cl, \VV}^{\infty}(M_0)$ are closed under taking adjoints.
\end{corollary}

We end this section with three remarks.

\begin{remark}\ Equation \eqref{eq.33} is
easily understood in the case of groups, when it amounts to the
possibility of covering any given compact set by finitely many
translations of a given open neighborhood of the identity. The
argument in general is the same as the argument used to define the
basic coordinate neighborhoods on $\GR$ in \cite{NistorINT}.  The
basic coordinate neighborhoods on $\GR$ were used in that paper to
define the smooth structure on the groupoid $\GR$.
\end{remark}

\begin{remark} We suspect that any proof of the fact that
$\Psi_{1,0,\VV}^{\infty}(M_{0})$ is closed under multiplication is
equivalent to the integrability of $A$. In fact, Melrose has
implicitly given some evidence for this in \cite{meicm} for particular
$(M, \VV)$, by showing that the kernels of the pseudodifferential
operators on $M_0$ that he constructed naturally live on a modified
product space $M^2_{\VV}$. In his case $M^2_{\VV}$ was a blow-up of
the product $M \times M$, and hence was a {\em larger}
compactification of the product $M_0 \times M_0$. The kernels of his
operators naturally extended to conormal distributions on this larger
product $M^2_{\VV}$. The product and adjoint were defined in terms of
suitable maps between $M^2_{\VV}$ and some fibered product spaces
$M^3_{\VV}$, which are suitable blow-ups of $M^3$ and hence larger
compactifications of $M_0^3$. This in principle leads to a solution of
the problem of microlocalizing $\VV$ that we stated in the
introduction whenever one can define the spaces $M^2_{\VV}$ and
$M^3_{\VV}$.  Let us also mention here that Melrose's approach usually
leads to algebras that are slightly larger than ours.
\end{remark}

\begin{remark} Let $\GR$ be a Lie groupoid 
such that the map $\piM$ is an isomorphism and let $N \subset M$ be
a face of $M$, then we obtain a generalized indicial map 
\begin{equation*}
  R_N : \Psi_{1, 0, \VV}^\infty(M_0) \simeq \Psi_{1, 0}^\infty(\GR)
  \to \Psi_{1, 0}^\infty(\GR_N).
\end{equation*}
In applications, the algebras $\Psi_{1, 0}^\infty(\GR_N)$ often turn out 
to be isomorphic to the
algebras $\Psi_{1, 0, \VV_1}^\infty(N_0; G)$ studied in the last
section of this paper. In fact, this is the motivation for introducing
the algebras $\Psi_{1, 0, \VV_1}^\infty(N_0; G)$.
\end{remark}

\section{Properties of $\Psi_{1,0, \VV}^\infty(M_0)$\label{sec.Prop}}
Theorem~\ref{theorem.gr} has several consequences similar to the
results in \cite{Mazzeo, MelroseScattering, meaps, MelroseMendoza,
Schrohe, schwil}.

\subsection{Basic properties}
We obtain that the algebras $\Psi_{1, 0, \VV}^\infty(M_{0})$ and
$\Psi_{cl, \VV}^\infty(M_{0})$ are independent of the choices made to
define them and, thus, depend only on the Lie structure at infinity
$(M, \VV)$.

\begin{corollary}\label{cor.indep}\
The spaces $\Psi_{1, 0, \VV}^m(M_{0})$ and $\Psi_{cl, \VV}^m(M_{0})$
are independent of the choice of the metric on $A$ and the function
$\chi$ used to define it, but depend, in general, on the Lie structure
at infinity $(M,A)$ on $M_0$.
\end{corollary}

\begin{proof}\
The space $\Psi_{1,0}^m(\GR)$ does not depend on the metric on $A$ or
on the function $\chi$ and neither does the vector representation
$\piM$. Then use Theorem \ref{theorem.gr}. The proof is the same
for classical operators.
\end{proof}

An important consequence is that $\Psi_{1, 0, \VV}^\infty(M_{0})$ 
and $\Psi_{cl, \VV}^\infty(M_{0}) = \cup_{m \in \ZZ}
\Psi_{cl, \VV}^m(M_{0})$ 
are filtered algebras, as it is the
case of the usual algebra of pseudodifferential operators on a compact
manifold.

\begin{proposition}\label{prop.alg}\
Using the above notation, we have that
\begin{multline*}
    \Psi_{1,0, \VV}^m(M_{0}) \Psi_{1, 0, \VV}^{m'}(M_0) \subseteq
    \Psi_{1, 0, \VV}^{m + m'}(M_{0}) \quad \text{and} \\ \Psi_{cl,
    \VV}^m(M_{0}) \Psi_{cl, \VV}^{m'}(M_0) \subseteq \Psi_{cl,
    \VV}^{m + m'}(M_{0})\,,
\end{multline*}
for all $m,m' \in \CC \cup \{-\infty\}.$
\end{proposition}

\begin{proof}\
Use Theorem \ref{theorem.gr} and the fact that $\piM$ preserves the
product.
\end{proof}

Part (i) of the following result is an analog of a standard result
about the $b$-calculus \cite{meaps}, whereas the second formula is the
independence of diffeomorphisms of the algebras $\Psi_{cl,
\VV}^\infty(M_{0})$, in the framework of manifolds with a Lie
structure at infinity. Recall that if $X\in\Gamma(A)$, we have denoted
by $\psi_{X} := \Psi_{X}(1,\cdot):M \rightarrow M$ the diffeomorphism
defined by integrating $X$ (and specializing at $t = 1$).

\begin{proposition} \label{prop.auto}\ (i) Let $x$ be a  defining
function of some hyperface of $M$. Then 
\begin{equation*}
x^{s} \Psi_{1,0,\VV}^m(M_{0}) x^{-s} = \Psi_{1, 0, \VV}^m(M_{0}) \quad \text{and} \quad 
x^{s} \Psi_{cl, \VV}^m(M_{0}) x^{-s} = \Psi_{cl, \VV}^m(M_{0})
\end{equation*} 
for any $s \in \CC$.

(ii) Similarly, 
$$
        \psi_X \Psi_{1, 0, \VV}^m(M_{0}) \psi_X^{-1} = \Psi_{1, 0,
        \VV}^m(M_{0})\quad \text{and} \quad \psi_X \Psi_{cl,
        \VV}^m(M_{0})\psi_X^{-1} = \Psi_{cl, \VV}^m(M_{0}),
$$ 
for any $X \in \Gamma(A)$.
\end{proposition}

\begin{proof}\
We have that $x^{s} \Psi^m_{cl}(\GR) x^{-s} = \Psi^m_{cl}(\GR)$, for
any $s \in \CC$, by \cite{LN1}. A similar result for type $(1,0)$
operators is proved in the same way as in \cite{LN1}. This proves (a)
because $\piM (x^s P x^{-s}) = x^s \piM(P) x^{-s}$.

Similarly, using the notations of Theorem \ref{theorem.gr}, we have
$\tilde{\psi}_X \Psi^m_{cl}(\GR) \tilde{\psi}_X^{-1} =
\Psi^m_{cl}(\GR)$, for any $X \in \Gamma(A) = \VV$. By the
diffeomorphism invariance of the space of pseudodifferential
operators, $\tilde{\psi}_X P \tilde{\psi}_X^{-1}$ defines a right
invariant family of pseudodifferential operators on $\GR$ for any such
right invariant family $P = (P_x)$, as in the proof of Theorem
\ref{theorem.gr}. To check that the family $P_1 := \psi_X P
\psi_X^{-1}$ has a compactly supported convolution kernel, denote by
$(\GR)_a = \{g, dist(g, d(g)) \le a\}$. Then observe that
$\supp(\tilde{\psi}_X P \tilde{\psi}_X^{-1}) \subset \GR_{d + 2\|X\|}$
whenever $\supp(P) \subset (\GR)_d$. Then use Equation
\eqref{eq.pi.psi} to conclude the result.

The proof for type $(1,0)$ operators is the same.
\end{proof}

Let us notice that the same proof gives (ii) above for any
diffeomorphism of $M_0$ that extends to an automorphism of $(M,A)$.
Recall that an automorphism of the Lie algebroid $\pi : A \to M$ is a
morphism of vector bundles $(\phi,\psi)$, $\phi : M \to M$, $\psi: A
\to A$, such that $\phi$ and $\psi$ are diffeomorphisms, $\pi \circ \psi =
\phi \circ \pi$, and we have
the following compatibility with the anchor map $\varrho$:
$$
    \varrho \circ \psi = \phi_* \circ \varrho.
$$

\subsection{Mapping properties}
Let $H^s(M_0)$ be the domain of $(1 + \Delta)^{s/2}$, where $\Delta$
is the (positive) Laplace operator on $M_0$ defined by the metric, if
$s \ge 0$. The space $H^{-s}(M_0)$, $s \ge 0$, is defined by duality,
the duality form being the pairing of distributions with test
functions.

\begin{corollary}\label{cor.bdd}\
Each operator $P \in \Psi^{m}_{1, 0, \VV}(M_{0})$, $P : \cunc(M_{0})
\rightarrow \cun(M_{0})$, extends to a continuous linear operators $P
: \cun(M) \rightarrow \cun(M)$ and $P: H^{s}(M_0) \to H^{s-m}(M_0)$.
The space $H^{m}(M_0)$, $m \ge 0$, identifies with the domain of $P$
with the graph topology and $H^{-m}(M_0) = PL^2(M_0) + L^2(M_0)$, for
any elliptic $P \in \Psi_{1, 0}^m(M_0)$.
\end{corollary}

\begin{proof}\
The first part is a direct consequence of the definitio of $\in
\Psi^{m}_{1, 0, \VV}(M_{0})$ is properly supported. The last part follows
  from the results of \cite{ain} and \cite{alnv}.

We now sketch the proof for the benefit of the reader. It follows from
the explicit form of the kernels of operators $T \in \Psi_{1, 0,
\VV}^{-n-1}(M_0)$, $n = \dim(M_0)$, that such a $T$ is bounded on
$L^2(M_0)$. Using the symbolic properties of the algebra $\Psi_{1, 0,
\VV}^\infty(M_0)$, namely Proposition \ref{prop.princ.symb} and
Proposition \ref{prop.alg}, it then follows that any $T \in \Psi_{1,
0, \VV}^{0}(M_0)$ is bounded on $L^2(M_0)$ (the details are the same
as in \cite{LN1} or \cite{alnv}). Using again the symbolic properties
of $\Psi_{1, 0, \VV}^\infty(M_0)$, it is proved as in \cite{alnv} that
the domain of the closure of $P$ and $PL^2(M_0) + L^2(M_0)$ are
independent of $P$ elliptic of order $m$.  Let us denote by $H_{m}$
the domain of the closure of $P$ and $H_{-m} = PL^2(M_0) +
L^2(M_0)$. Then it is proved in \cite{alnv} that $T : H_{r} \to
H_{r-m}$ is bounded, for any $T$ of order $m$. In \cite{ain} it is
proved using partitions of unity that $T : H^{r}(M_0) \to
H^{r-m}(M_0)$ is bounded for any $T$ of order $m$. This shows that
$H_r= H^r(M_0)$ for any $r \in \RR$.
\end{proof}

\subsection{Quantization} We have the following quantization
properties of the algebra $\Psi_{1, 0, \VV}^\infty(M_0)$.

For any $X \in \Gamma(A)$, denote by $a_X : A^* \to \CC$ the function
defined by $a_X(\xi) = \xi(X)$. Then there exists a unique Poisson
structure on $A^*$ such that $\{a_X, a_Y\} = a_{[X, Y]}$.  It is
related to the Poisson structure $\{\,\cdot\,,\,\cdot\,\}^{T^*M}$ on
$T^*M$ via the formula
$$
        \{f_1\circ \rho^*,f_2\circ \rho^*\}^{T^*M}
        = \{f_1,f_2\}\circ \rho^*,
$$
where $\rho^*:T^*M\to A^*$ denotes the dual to the anchor map $\rho$.
In particular, in $\{\,\cdot\,,\,\cdot\,\}$ and
$\{\,\cdot\,,\,\cdot\,\}^{T^*M}$ coincide on $M_0$.

\begin{proposition}\label{prop.comm}\
We have that
\begin{equation*}
    \sigma^{(m + m' - 1)}([P, Q]) = \{ \sigma^{(m)}(P),
    \sigma^{(m')}(Q) \}
\end{equation*}
for any $P \in \Psi^{m}_{1,0}(M_0)$ and any $Q \in
\Psi^{m'}_{1,0}(M_0)$, where $\{\,\cdot \,,\, \cdot \,\}$ is the usual
Poisson bracket on $A^*$.
\end{proposition}

\begin{proof}\ The Poisson structure on $T^*M_0$ is induced from 
the Poisson structure on $A^*$. In turn, the Poisson structure on
$T^*M_0$ determines the Poisson structure on $A^*$, because $T^*M_0$ 
is dense in $A^*$.  The desired result then follows from the
similar result that is known for pseudodifferential operators on $M_0$
and the Poisson bracket on $T^*M_0$.
\end{proof}

We conclude with the following result, which is independent of the
previous considerations, but sheds some light on them. The invariant
differential operators on $\GR$ are generated by $d$--vertical
invariant vector fields on $\GR$, that is by $\Gamma(A(\GR))$. We have
by definition that $\piM = \varrho : \Gamma(M; A(\GR)) \to \Gamma(M;
TM)$, and hence $\piM$ maps the algebra of invariant differential
operators onto $\GR$ to $\DiffV{*}(M_0)$. In particular, the proof of
Theorem \ref{theorem.gr} (more precisely Equation \eqref{eq.surj}) can
be used to prove the following result, which we will however prove
also without making appeal to Theorem \ref{theorem.gr}.

\begin{proposition}\label{prop.quant}\
Let $X \in \Gamma(A)$ and denote by $a_{X} (\xi) = \xi(X)$ the
associated linear function on $A^*$. Then $a_{X,\chi} \in S^1(A^*)$ and
$a_X(D) = -i  X$. Moreover,
$$ 
        \{a_\chi(D), \, a =\mbox{ \rm polynomial in each fiber}\,\} =
        \DiffV{*}(M_0).
$$
\end{proposition}

\begin{proof}\
We continue to use a fixed metric on $A$ to trivialize any density
bundle. Let $u = \ffi(a)$, where $a \in S^m_{cl}(A^*)$ polynomial in
each fiber. By the Fourier inversion formula (and integration by
parts), $u$ is supported on $M$, which is the same thing as saying
that $u$ is a distribution the form $\< u, f \> = \int_{M} P_0 f(x)
d\vol(x)$, with $P_0$ a differential operator acting along the fibers
of $A \to M$ and $f \in \CIc(A)$. It then follows from the definition
of $a_\chi(D)$, from the formula above for $u = \ffi(a)$, and from the
fact that $\chi = 1$ in a neighborhood of the support of $u$ that
\begin{equation}\label{eq.form.aD}
    a_\chi(D)f(x) = [P_0f(\exp_x(-v))] \vert_{v = 0}, \quad v \in
    T_xM_0.
\end{equation}

Let $X_1,X_2, \ldots, X_m \in \Gamma(A)$ and
\begin{equation}\label{eq.a.prod}
    a = a_{X_1} a_{X_2} \ldots a_{X_m} \in S^m(A^*).
\end{equation}
Then the differential operator $P_0$ above is given by the formula
\begin{equation*}
    P_0f(x) = \int_{A_x^*} a(\xi) {\mathcal F}^{-1}f(\xi),
\end{equation*}
with the inverse Fourier transform ${\mathcal F}^{-1}$ being defined
along the fiber $A_x$. Hence
\begin{equation*}
    P_0 = i ^{m} X_1 X_2 \ldots X_m,
\end{equation*}
with each $X_j$ being identified with the family of constant
coefficients differential operators along the fibers of $A \to M$ that
acts along $A_x$ as the derivation in the direction of $X_j(x)$.

For any $X \in A$, we shall denote by $\psi_{tX}$ the one parameter
subgroup of diffeomorphisms of $M$ generated by $X$.  (Note that
$\psi_{tX}$ is defined for any $t$ because $M$ is compact and $X$ is
tangent to all faces of $M$.) We thus obtain an action of $\psi_{tX}$
on functions by $[\psi_{tX}(f)](x) = f(\exp(tX)x)$. Then the
differential operator $P_0$ is associated to $a$ as in Equation
\eqref{eq.a.prod} is given by
\begin{equation}\label{eq.above}
    P_0(f \circ \exp)\vert_M = i ^{m} \big[ \pa_1 \pa_2 \ldots
    \pa_m \psi_{t_1X_1 + t_2 X_2 + \ldots + t_m X_m} f \big]\vert_{t_1
    = \ldots = t_m = 0}.
\end{equation}

Then Equations \eqref{eq.form.aD} and \eqref{eq.above} give
\begin{equation}
    a_\chi(D)f =i ^{m}\big[ \pa_1 \ldots \pa_m \exp (-t_1X_1 -
    \ldots - t_m X_m) f \big]\vert_{t_1 = \ldots = t_m = 0},
\end{equation}
In particular, $a_X(D) = -i  X$, for any $X \in \Gamma(A)$.

This proves that
\begin{equation}\label{eq.ainD}
    a_\chi(D) \in \DiffV{*}(M_0),
\end{equation}
by the Campbell-Hausdorff formula \cite{hochschild, onishchik}, which
states that $a_\chi(D)$ is generated by $X_1, X_2, \ldots, X_n$ (and
their Lie brackets), and hence that it is generated by $\VV$, which
was assumed to be a Lie algebra.

Let us prove now that any differential operator $P \in \DiffV{*}(M_0)$
is of the form $a_\chi(D)$, for some polynomial symbol $a$ on
$A^*$. This is true if $P$ has degree zero. Indeed, assume $P$ is the
multiplication by $f \in \CI(M)$. Lift $f$ to an order zero symbol on
$A^*$, by letting this extension to be constant in each fiber. Then $P
= f(D)$. We shall prove our statement by induction on the degree $m$
of $P$. By linearity, we can reduce to the case $P = i ^{-m}X_1
\ldots X_m$, where $X_1, \ldots, X_m \in \Gamma(A)$. Let $a = a_{X_1}
\ldots a_{X_m}$. Then
\begin{equation*}
        \sigma_m(a_\chi(D))(\xi) = a(\xi) =
        X_1(\xi) \ldots X_m(\xi) = \sigma_m(P),
\end{equation*}
and hence $Q := a_\chi(D) - i^{-m} X_1 \ldots X_m \in
\DiffV{m-1}(M_0)$.  By the induction hypothesis, $Q = b_{\chi}(D)$ for
some polynomial symbol of order at most $m-1$ on $A^*$. This completes
the proof.
\end{proof}

{}From this we obtain the following corollary.

\begin{corollary}\label{cor.quant}\
Let $\Diff{M_0}$ be the algebra of all differential operators on
$M_0$. Then
\begin{equation*}
    \Psi_{1,0,\VV}^\infty(M_0) \cap \Diff{M_0} = \DiffV{*}(M_0).
\end{equation*}
\end{corollary}

\begin{proof}\ We know from the above proposition that
$$
        \Psi_{1,0,\VV}^\infty(M_0) \cap \Diff{M_0} \supset
        \DiffV{*}(M_0).
$$ 
Conversely, assume $P \in \Psi_{1,0,\VV}^m(M_0) \cap \Diff{M_0}$. We
shall prove by induction on $m$ that $P \in \DiffV{m}(M_0)$. If $m =
0$ then $P$ is the multiplication with a smooth function $f$ on
$M_0$. But then $f = \sigma^{(0)}(P) \in S^0(A^*)$ is constant along
the fibers of $A^* \to M$, and hence $f \in \CI(M) \subset
\DiffV{*}(M_0)$. 

Assume now that the statement is proved for $P$ of
order $< m$. We shall prove it then for $P$ of order $m$. Then $a :=
\sigma^{(m)}(P)$ is a polynomial symbol in $S^m(A^*)$. Thus $a_\chi(D)
\in \DiffV{m}(M_0)$, by Proposition \ref{prop.quant}. But then
$\sigma^{(m)}(P - a_\chi(D)) = 0$, by Lemma \ref{basqu}, and hence $P
- a_\chi(D) \in \Psi_{1,0,\VV}^{m-1}(M_0) \cap \Diff{M_0}$. By the
induction hypothesis $P - a_\chi(D) \in \DiffV{m-1}(M_0)$. This
completes the proof.
\end{proof}


\section{Group actions and semi-classical limits\label{Sec.Ginv}}

One of the most convenient features of manifolds with a Lie structure
at infinity is that questions on the analysis on these manifolds often
reduce to questions on the analysis on simpler manifolds. These
simpler manifolds are manifolds of the same dimension but endowed with
certain non-trivial group actions.  Harmonic analysis techniques then
allow us to ultimately reduce our questions to analysis on lower
dimensional manifolds with a Lie structure at infinity. In this
section, we discuss the algebras $\Psi_{1,0,\VV}^\infty(M_0, G)$ that
generalize the algebras $\Psi_{1,0,\VV}^\infty(M_0)$ when group
actions are considered.  These algebras are necessary for the
reductions mentioned above and are typically the range of
(generalized) indicial maps.  Then we discuss a semi-classical version
of the algebra $\Psi_{1,0,\VV}^\infty(M_0)$.

\subsection{Group actions} We shall consider
the following setting. Let $M_0$ be a manifold with a Lie
structure at infinity $(M,A)$, and $\VV = \Gamma(A)$, as above.
Also, let $G$ be a Lie group with Lie algebra $\fgg := \Lie(G)$.
We shall denote by $\fgg_M$ the bundle $M \times \fgg \to M$. Then
\begin{equation}
    \VV_G : = \VV \oplus \CI(M, \fgg) \simeq
    \Gamma(A \oplus \fgg_M)
\end{equation}
has the structure of a Lie algebra with respect to the bracket
$[\,\cdot \,,\, \cdot \,]$ which is defined such that on $\CI(M,
\fgg)$ it coincides with the pointwise bracket, on $\VV$ it
coincides with the original bracket, and, for any $X \in \VV$, $f
\in \CI(M)$, and $Y \in \fgg$, we have
\begin{equation*}
    [X, f \otimes Y] := X(f) \otimes Y\,.
\end{equation*}
(Here $f \otimes Y$ denotes the function $\xi : M \to \fgg$
defined by $\xi(m) = f(m)Y \in \fgg$.)

The main goal of this subsection is to indicate how the results of
the Section \eqref{sec.KN} extend to $\VV_G$, after we replace $A$
with $A \oplus \fgg_M$, $M_0$ with $M_0 \times G$, and $M$ with $M
\times G$. The resulting constructions and definitions will yield
objects on $M \times G$ that are invariant with respect to the action of $G$ on
itself by {\em right} translations.

We now proceed by analogy with the construction of the operators
$a_{\chi}(D)$ in Subsection \ref{subsec.KN}. First, we identify a
section of
$\VV_G := \VV \oplus \CI(M, \fgg) \simeq \Gamma(A \oplus \fgg_M)$
with a right $G$-invariant vector field on $M_0 \times G$. At the
level of vector bundles, this corresponds to the map
\begin{equation}\label{eq.def.p}
    p: T(M_0 \times G) = TM_0 \times TG \to TM_0\times \fgg,
\end{equation}
where the map $TG \to \fgg$ is defined using the trivialization of
$TG$ by {\em right} invariant vector fields. Let $p_1 : M \times G
\to M$ be the projection onto the first component and $p_1^*A$ be
the lift of $A$ to $M \times G$ via $p_1$.

The map $p$ defined in the Equation \eqref{eq.def.p} can then be
used to define the lift
\begin{equation}\label{eq.p*}
    p^*(u) \in I^m(p_1^*A \oplus TG, M \times G),
\end{equation}
for any distribution $u \in I^m(A \oplus \fgg_M, M)$. In
particular, $p^*(u)$ will be a right $G$-invariant distribution.
Then we define $\mR$ to be the restriction of distributions from
$p_1^*A \oplus TG$ to distributions on $TM_0 \times TG = T(M_0
\times G)$.

We endow $M_0 \times G$ with the metric obtained from a metric on
$A$ and a right invariant metric on $G$. This allows us to define
the exponential map, thus obtaining, as in Section \ref{sec.KN}, a
differentiable map
\begin{equation}\label{eq.def.PhiG}
    \Phi : (TM_0 \times TG)_r = (T(M_0 \times G))_r \to (M_0 \times G)^2
\end{equation}
that is a diffeomorphism onto an open neighborhood of the
diagonal, provided that $r < r_0$, where $r_0$ is the injectivity
radius of $M_0 \times G$. We shall denote as before by
$$
    \Phi_* : I^m_c( (TM_0 \times TG)_r  , M_0 \times G) \to
    I^m_c((M_0 \times G)^2, M_0 \times G)
$$
the induced map on conormal distributions.

The inverse Fourier transform will give a map
\begin{equation}\label{eq.FinvG}
    \cF_{{\rm fiber}}^{-1}: S^{m}_{1,0}(A^{*} \oplus \fgg_M^*)
    \longrightarrow I^m(A \oplus \fgg_M, M),
\end{equation}
defined by the same formula as before (Equation \eqref{eq.Finv}).
Finally, we shall also need a smooth function $\chi$ on $A \oplus
\fgg_M$ that is equal to 1 in a neighborhood of the zero section
and has support inside $(A \oplus \fgg_M)_r$.

We can then define the quantization map in the $G$-equivariant
case by
\begin{equation}\label{defquant.G}
    q_{\Phi,\chi,G} := \Phi_{*} \circ \cR \circ p^* \circ \chi
    \circ \cF_{{\rm fiber}}^{-1} : S^{m}_{1,0} (A^{*} \oplus
    \fgg_M^*)\longrightarrow \mathcal I^m( (M_0 \times G)^2, M_0
    \times G).
\end{equation}
The main difference with the definition in Equation
\eqref{defquant} is that we included the map $p^*$, which is the
lift of distributions in $I^m(A \oplus \fgg_M, M)$ to
$G$-invariant distributions in $I^m(p_1^*A \oplus TG, M_0 \otimes
G)$, see Equation \eqref{eq.p*}. Then
\begin{equation}\label{eq.def.aDG}
    a_\chi(D) = T \circ q_{\Phi, \chi, G},
\end{equation}
as before.

With this definition of the quantization map, all the results of
the previous sections remain valid, with the appropriate
modifications. In particular, we obtain the following definition
of the algebra of $G$-equivariant pseudodifferential operators
associated to $(A , M, G)$.

\begin{definition}\label{def.psi.MAG}\
For $m\in\RR$, the space $\Psi_{1,0,\VV}^{m}(M_{0}, G)$ of {\em
$G$-equivariant pseudodifferential operators generated by the Lie
structure at infinity $(M,A)$} is the linear space of operators
$\cunc(M_0 \times G) \rightarrow \cunc(M_{0} \times G)$ generated
by $a_\chi(D)$, $a \in S_{1,0}^m(A^* \oplus \fgg_M^*)$, and
$b_{\chi}(D)\psi_{X_1}\ldots \psi_{X_k}$, $b \in S^{-\infty}(A^*
\oplus \fgg_M^*)$ and $X_j \in \Gamma(A \oplus \fgg_M)$.

The space $\Psi_{cl, \VV}^{m}(M_{0}, G)$ of {\em classical
$G$-equivariant pseudodifferential operators generated by the Lie
structure at infinity $(M,A)$} is defined similarly, but using
classical symbols $a$.
\end{definition}

With this definition, all the results on the algebras
$\Psi_{1,0,\VV}^{m}(M_{0})$ and $\Psi_{cl,\VV}^{m}(M_{0})$ extend
right away to the spaces $\Psi_{1,0,\VV}^{m}(M_{0}, G)$ and
$\Psi_{cl,\VV}^{m}(M_{0}, G)$. In particular, these spaces are
algebras, are independent of the choice of the metric on $A$ used
to define them, and have the usual symbolic properties of the
algebras of pseudodifferential operators.

The only thing that maybe needs more explanations is with what we
replace $\piM$ in the $G$-equivariant case, because in the
$G$-equivariant case we no longer use the vector representation.
Let $\GR$ be a groupoid integrating $A$, $\Gamma(A) = \VV$. Then $\GR \times G$
integrates $A \oplus \fgg_M$.  If $P = (P_x) \in \Psi^m_{1,0} (\GR
\times G)$, then we consider $\pi_0(P)$ to be the operator induced
by $P_x$ on $(\GR_x/\GR_x^x) \times G$, $x \in M_0$, the later
space being a quotient of $(\GR \times G)_x$. We shall use then
$\pi_0$ instead of $\piM$ in the $G$-equivariant case. (By the
proof of Theorem \ref{theorem.gr}, $\pi_0 = \piM$, if $G$ is
reduced to a point.)

\subsection{Indicial maps}
The main reason for considering the algebras
$\Psi_{1,0,\VV}^{m}(M_{0}, G)$ and their classical counterparts is the
following. Let $(M, \VV)$, $\VV = \Gamma(M, A)$, be a manifold with a
Lie structure at infinity. Let $N_0 \subset M$ be a submanifold such
that $T_x N_0 = \varrho(A_x)$ for any $x \in N_0$. Moreover, assume
that $N_0$ is completely contained in an open face $F \subset M$ such
that $N := \overline{N_0}$ is a submanifold with corners of $F$ and
$N_0 = N \smallsetminus \pa N$.  Then the restriction $A\vert_{N_0}$
is such that the Lie bracket on $\VV = \Gamma(A)$ descends to a Lie
bracket on $\Gamma(A\vert_{N_0})$. (This is due to the fact that the
space $I$ of functions vanishing on $N$ is invariant for derivations
in $\VV$. Then $I\VV$ is an ideal of $\VV$, and hence $\VV/I\VV \simeq
\Gamma(A\vert_{N})$ is naturally a Lie algebra.)

Assume now that there exists a Lie group $G$ and a vector bundle $A_1
\to N$ such that $A \vert_{N} \simeq A_1 \oplus \fgg_N$ and $\VV_1 :=
\VV\vert_{N} \simeq \Gamma(A_1)$.  Then $\VV_1$ is a Lie algebra and
$(N_0, N, A_1)$ is also a manifold with a Lie structure at infinity.
In many cases (certainly for many of the most interesting examples)
one obtains for any Lie group $H$ a natural morphism
\begin{equation}\label{eq.def.RN}
    R_N : \Psi_{1, 0, \VV}^{m}(M_0; H) \to \Psi_{1, 0,
    \VV_1}^{m}(N_0; G \times H).
\end{equation}
For example, the generalizations of the morphisms considered in
\cite{LN1} are of the form~\eqref{eq.def.RN}.  However, we do not know
exactly what are the conditions under which the morphism $R_N$ above
is defined.

Let $\mathfrak h = \Lie H$ and $\mathfrak h_N = M \times \Lie H$.
Then, at the level of kernels the morphism defined by Equation
\eqref{eq.def.RN} corresponds to the restriction maps
\begin{equation*}
    r_N : I^m(A^* \oplus \mathfrak h_N^*, M) \to I^*(A^*\vert_{N}
    \oplus \mathfrak h_N^*, N) \simeq I^*(A_1^* \oplus \fgg_N
    \oplus \mathfrak h_N^*, N)
\end{equation*}
in the sense that $R_N(a_\chi(D)) = (r_N(a))_\chi(D)$.

\subsection{Semi-classical limits} We now define the algebra
$\Psi_{1,0,\VV}^{m}(M_{0}[[h]])$, an element of which will be, roughly
speaking, a semi-classical family of operators $(T_t)$, $T_t \in
\Psi_{1,0,\VV}^{m}(M_{0})$ $t \in (0, 1]$.
See \cite{VZ} for some applications of semi-classical analysis.

\begin{definition}\label{def.psi.MAAD}\
For $m\in\RR$, the space $\Psi_{1,0,\VV}^{m}(M_{0}[[h]])$ of {\em
pseudodifferential operators generated by the Lie structure at
infinity $(M,A)$} is the linear space of families of operators $T_t:
\cunc(M_0 \times G) \rightarrow \cunc(M_{0} \times G)$, $t \in (0,
1]$, generated by
\begin{equation*}
    a_\chi(t, t D)\,, \quad a \in S_{1,0}^m([0, 1] \times A^* \oplus
    \fgg_M^*),
\end{equation*}
and
\begin{equation*}
    b_{\chi}(t, t D)\psi_{tX_1(t)}\ldots \psi_{tX_k(t)}, \quad b \in
    S^{-\infty}([0, 1] \times A^* \oplus \fgg_M^*), \quad X_j \in
    \Gamma([0, 1] \times A \oplus \fgg_M).
\end{equation*}
\end{definition}

The space $\Psi_{cl, \VV}^{m}(M_{0}[[h]])$ of {\em semi-classical
families of pseudodifferential operators generated by the Lie
structure at infinity $(M,A)$} is defined similarly, but using
classical symbols $a$.

Thus we consider families of operators $(T_t)$, $T_t \in
\Psi_{1,0,\VV}^{m}(M_{0})$, defined in terms of data $a, b, X_k$, that
extends smoothly to $t = 0$, with the interesting additional feature
that the cotangent variable is rescaled as $t \to 0$.

Again, all the results on the algebras $\Psi_{1,0,\VV}^{m}(M_{0})$ and
$\Psi_{cl,\VV}^{m}(M_{0})$ extend right away to the spaces
$\Psi_{1,0,\VV}^{m}(M_{0}[[h]])$ and $\Psi_{cl,\VV}^{m}(M_{0}[[h]])$,
except maybe Proposition \ref{prop.quant} and its corollary, Corollary
\ref{cor.quant}, that need to be properly reformulated.

Another variant of the above constructions is to consider families of
manifolds with a Lie structure at infinity. The necessary changes are
obvious though, and we will not discuss them here.


\begin{thebibliography}{10}

\bibitem{ain} B.~Ammann, A.~Ionescu, and V.~Nistor, \emph{Sobolev
spaces and regularity for polyhedral domains}, Preprint 2004/2005,
ArXiv:math.AP/0402321.

\bibitem{aln1} B.~Ammann, R.~Lauter, and V.~Nistor, \emph{On the
geometry of {R}iemannian manifolds with a {L}ie structure at
infinity}, Int. J. Math.  Math. Sci.  (2004), no.~1-4, 161--193.

\bibitem{alnv} B.~Ammann, R.~Lauter, V.~Nistor, and A.~Vasy,
\emph{Complex powers and non-compact manifolds},
Comm. Part. Diff. Eq. \textbf{29} (2004), no.~5-6, 671--705.

\bibitem{BNXZ} C.~Bacuta, V.~Nistor, and L.~Zikatanov, \emph{Improving
the rate of convergence of `high order finite elements' on polygons},
Numerische Mathematik \textbf{100} (2005), no.~2, 165--184.

\bibitem{Baer} C.~B{\"a}r, \emph{The {D}irac operator on hyperbolic
manifolds of finite volume}, J. Diff. Geom. \textbf{54} (2000), no.~3,
439--488.

\bibitem{Costabel} M.~Costabel, \emph{Boundary integral operators on
curved polygons}, Ann. Mat.  Pura Appl. (4) \textbf{133} (1983),
305--326.

\bibitem{CrainicFernandes} M.~Crainic and R.~Fernandes,
\emph{Integrability of {L}ie brackets}, Ann. of Math. (2) \textbf{157}
(2003), no.~2, 575--620.

\bibitem{Dauge} M.~Dauge, \emph{Elliptic boundary value problems on
corner domains}, Lecture Notes in Mathematics, vol. 1341,
Springer-Verlag, Berlin, 1988, Smoothness and asymptotics of
solutions.

\bibitem{emm91} C.~Epstein, R.B. Melrose, and G.~Mendoza,
\emph{Resolvent of the {L}aplacian on strictly pseudoconvex domains},
Acta Math. \textbf{167} (1991), 1--106.

\bibitem{hochschild} G.~Hochschild, \emph{{The structure of Lie
groups.  (Holden-Day Series in Mathematics)}}, {San
Francisco-London-Amsterdam: Holde-Day, Inc. IX, 230 p.  }, 1965.

\bibitem{fio} L.~H\"ormander, \emph{Fourier integral operators {I}},
Acta Math. \textbf{127} (1971), 79--183.

\bibitem{hor3} L.~H\"ormander, \emph{The analysis of linear partial
differential operators, {\rm vol.\ 3}. {P}seudo-differential
operators}, Grundlehren der Mathematischen Wissenschaften, vol. 274,
Springer-Verlag, Berlin - Heidelberg - New York, 1985.

\bibitem{Karoubi} M.~Karoubi, \emph{Homologie cyclique et
{K}-theorie}, Ast\'erisque \textbf{149} (1987), 1--147.

\bibitem{Spencer} A.~Kumpera and D.~Spencer, \emph{{Lie Equations}},
vol.~{1}, {Princeton Univ.  Press.}, 1972.

\bibitem{zfr} R.~Lauter, \emph{Pseudodifferential analysis on
conformally compact spaces}, Mem.\ Amer.\ Math.\ Soc., vol. 163, 2003.


\bibitem{defr} R.~Lauter and S.~Moroianu, \emph{Fredholm theory for
degenerate pseudodifferential operators on manifolds with fibered
boundaries}, Comm.  Part. Diff. Eq. \textbf{26} (2001), 233--283.

\bibitem{LN1} R.~Lauter and V.~Nistor, \emph{Analysis of geometric
operators on open manifolds: a groupoid approach}, Quantization of
Singular Symplectic Quotients (N.P. Landsman, M.~Pflaum, and
M.~Schlichenmaier, eds.), Progress in Mathematics, vol. 198,
Birkh{\"a}user, Basel - Boston - Berlin, 2001, pp.~181--229.

\bibitem{Lewis} J.~Lewis, \emph{Layer potentials for elastostatics and
hydrostatics in curvilinear polygonal domains},
Trans. Amer. Math. Soc. \textbf{320} (1990), no.~1, 53--76.

\bibitem{LewisParenti} J.~Lewis and C.~Parenti,
\emph{Pseudodifferential operators of {M}ellin type},
Comm. Part. Diff. Eq. \textbf{8} (1983), no.~5, 477--544.

\bibitem{Mackenzie} K.~Mackenzie, \emph{{Lie groupoids and {L}ie
algebroids in differential geometry}}, {Lecture Notes Series},
vol. {124}, {London Mathematical Society}, 1987.

\bibitem{Mazzeo} R.~Mazzeo, \emph{{Elliptic theory of differential
edge operators. I.}}, Comm.  Part. Diff. Eq. \textbf{16} (1991),
no.~10, 1615--1664.

\bibitem{mame87} R.~Mazzeo and R.B. Melrose, \emph{Meromorphic
extension of the resolvent on complete spaces with asymptotically
constant negative curvature}, J. Funct.  Anal. \textbf{75} (1987),
260--310.

\bibitem{MaMeAsian} R.~Mazzeo and R.B. Melrose,
\emph{Pseudodifferential operators on manifolds with fibered
boundaries}, Asian J. Math. \textbf{2} (1998), 833--866.

\bibitem{MazzeoVasy} R.~Mazzeo and A.~Vasy, \emph{Analytic
continuation of the resolvent of the {L}aplacian on {$SL(3)/SO(3)$}},
Amer. J. Math. \textbf{126} (2004), 821--844.
 
\bibitem{MelroseCorners} R.B. Melrose, \emph{Analysis on manifolds
with corners}, book: in preparation.

\bibitem{me81} R.B. Melrose, \emph{Transformation of boundary value
problems}, Acta Math.  \textbf{147} (1981), 149--236.

\bibitem{meicm} R.B. Melrose, \emph{Pseudodifferential operators,
corners and singular limits}, Proceeding of the International Congress
of Mathematicians, Kyoto (Berlin - Heidelberg - New York),
Springer-Verlag, 1990, pp.~217--234.

\bibitem{meaps} R.B. Melrose, \emph{{The Atiyah-Patodi-Singer index
theorem.}}, {Research Notes in Mathematics (Boston,
Mass.). 4. Wellesley, MA: A. K. Peters, Ltd.}, 1993.

\bibitem{MelroseScattering} R.B. Melrose, \emph{Geometric scattering
theory}, Stanford Lectures, Cambridge University Press, Cambridge,
1995.

\bibitem{MelroseMendoza} R.B. Melrose and G.~Mendoza, \emph{Elliptic
operators of totally characteristic type}, MSRI Preprint.

\bibitem{mendoza} G. Mendoza, \emph{Strictly pseudoconvex {$b$}-{CR}
manifolds}, Comm.  Part. Diff. Eq. \textbf{29} (2004), no.~9-10,
1437--1503.

\bibitem{MiNi} M.~Mitrea and V.~Nistor, \emph{Boundary layer
potentials on manifolds with cylindrical ends}, ESI Preprint no. 1244,
2002.

\bibitem{NagelStein} A.~Nagel and E.~Stein, \emph{Lectures on
pseudodifferential operators: regularity theorems and applications to
nonelliptic problems}, Mathematical Notes, vol.~24, Princeton
University Press, Princeton, N.J., 1979.

\bibitem{NistorINT} V.~Nistor, \emph{Groupoids and the integration of
{L}ie algebroids}, J. Math.  Soc. Japan \textbf{52} (2000), 847--868.

\bibitem{nwx} V.~Nistor, A.~Weinstein, and P.~Xu,
\emph{Pseudodifferential operators on groupoids}, Pacific
J. Math. \textbf{189} (1999), 117--152.

\bibitem{onishchik} A.~Onishchik, \emph{{Lie groups and Lie algebras
I. Foundations of Lie theory. Lie transformation groups. Transl. from
the Russian by A.  Kozlowski.}}, {Encyclopaedia of Mathematical
Sciences. 20. Berlin: Springer-Verlag.}, 1993.

\bibitem{Parenti} C.~Parenti, \emph{Operatori pseudodifferentiali in
{$\RR^n$} e applicazioni}, Annali Mat. Pura ed App. \textbf{93}
(1972), 391--406.

\bibitem{Schrohe} E.~Schrohe, \emph{{Spectral invariance, ellipticity,
and the Fredholm property for pseudodifferential operators on weighted
Sobolev spaces.}}, Ann. Global Anal. Geom. \textbf{10} (1992), no.~3,
237--254.

\bibitem{ScSc} E.~Schrohe and B.W. Schulze, \emph{Boundary value
problems in {B}outet de {M}onvel's algebra for manifolds with conical
singularities {II}}, Boundary value problems, Schr\"odinger operators,
deformation quantization, Math.  Top., vol.~8, Akademie Verlag,
Berlin, 1995.

\bibitem{schwil} B.W. Schulze, \emph{{Boundary value problems and
singular pseudo-differential operators.}}, {Wiley-Interscience Series
in Pure and Applied Mathematics.  Chichester: John Wiley \& Sons.},
1998.

\bibitem{ShubinBook} M. Shubin, \emph{{Pseudodifferential operators
and spectral theory}}, {Springer Verlag}, Princeton, N.J., 1987.

\bibitem{sim} S. Simanca, \emph{Pseudo-differential operators}, Pitman
research notes in mathematics, vol. 171, Longman Scientific \&
Technical, Harlow, Essex, 1990.

\bibitem{Taylor1} M.~Taylor, \emph{Pseudodifferential operators},
Princeton Mathematical Series, vol.~34, Princeton University Press,
Princeton, N.J., 1981.

\bibitem{Taylor2} M.~Taylor, \emph{Partial differential equations},
Applied Mathematical Sciences, vol. I--III, Springer-Verlag, New York,
1995-1997.

\bibitem{VasyN} A.~Vasy, \emph{Propagation of singularities in
many-body scattering}, Ann. Sci.  \'Ecole Norm. Sup. (4) \textbf{34}
(2001), no.~3, 313--402.

\bibitem{VZ} A.~Vasy and M.~Zworski, \emph{Semiclassical estimates in
asymptotically {E}uclidean scattering}, Comm. Math. Phys. \textbf{212}
(2000), no.~1, 205--217.

\bibitem{jaredduke} J.~Wunsch, \emph{Propagation of singularities and
growth for {S}chr{{\"o}}dinger operators}, Duke Math.~J. \textbf{98}
(1999), 137--186.
\end{thebibliography}
\end{document}